\documentclass{article}

\def\mathbb{\Bbb}
\usepackage{amsfonts,latexsym,amstext}
\usepackage{amsmath}
\usepackage{amssymb}

\usepackage{comment}
\usepackage{amsthm}
\usepackage[latin1]{inputenc}
\usepackage{color}

\newcommand{\abs}[1]{\left| #1 \right|}                                           

\theoremstyle{plain}
\newtheorem{theorem}{Theorem}[section]

\newtheorem{lemma}[theorem]{Lemma}
\newtheorem{proposition}[theorem]{Proposition}
\newtheorem{corollary}[theorem]{Corollary}

\theoremstyle{definition}
\newtheorem{definition}[theorem]{Definition}

\newtheorem{example}[theorem]{Example}
\newtheorem{hypothesis}[theorem]{Hypothesis}

\theoremstyle{remark}
\newtheorem{remark}[theorem]{Remark}

\makeatletter
\@addtoreset{equation}{section}

\makeatother

\def\sqr#1#2{{\vcenter{\vbox{\hrule height .#2pt \hbox{\vrule
 width .#2pt height#1pt \kern#1pt \vrule
width .#2pt} \hrule height .#2pt}}}}

\def\ds{\begin{displaystyle}}
\def\eds{\end{displaystyle}}
\def\dis{\displaystyle }
\def\<{\langle }
\def\>{\rangle }

\def\R{\mathbb R}

\def\E{\mathbb E}
\def\P{\mathbb P}

\def\cala{{\cal A}}
\def\calb{{\cal B}}

\def\calf{{\cal F}}

\def\calk{{\cal K}}

\newcommand{\probxa}[1]{\mathbb{P}^{x,a} \left\{ #1 \right\}}                           
\newcommand{\sperxa}[1]{\mathbb{E}^{x,a} \left[ #1 \right]}                               
\newcommand{\sperxaut}[1]{\mathbb{E}_{u,t}^{x,a} \left[ #1 \right]}               

\DeclareMathAlphabet{\mathonebb}{U}{bbold}{m}{n}                           %
\newcommand{\one}{\ensuremath{\mathonebb{1}}}                               

\title{Optimal control of semi-Markov processes with a backward stochastic differential equations approach}

\date{}
\author{Elena Bandini\thanks{Politecnico di Milano,
Dipartimento di Matematica, piazza Leonardo da Vinci 32, 20133 Milano, Italy.
e-mail: elena.bandini@polimi.it.}
\, and \, Fulvia Confortola\thanks{Politecnico di Milano,
Dipartimento di Matematica, piazza Leonardo da Vinci 32, 20133 Milano, Italy.
e-mail: fulvia.confortola@polimi.it.}}

\begin{document}
\allowdisplaybreaks

\maketitle

\begin{abstract}
\noindent In the present work we employ, for the first time, backward stochastic differential equations (BSDEs) to study the optimal control of semi-Markov processes on finite horizon, with general state and action spaces. More precisely, we prove that the value function and the optimal control law can be represented by means of the solution of a class of BSDEs driven by a semi-Markov  process  or,
equivalently, by the associated random measure.
The peculiarity of the semi-Markov framework, with respect to the pure jump Markov case, consists in the proof of the relation between BSDE and optimal control problem. This is done, as usual, via the Hamilton-Jacobi-Bellman (HJB) equation, which however in the semi-Markov case is characterized by an additional differential term $\partial_a$. Taking into account the particular structure of semi-Markov processes we rewrite the HJB equation in a suitable integral form which involves a directional derivative operator $D$ related to $\partial_a$. Then, using a formula of It$\hat{\mbox{o}}$ type tailor-made for semi-Markov processes and the operator $D$, we are able to prove that the BSDE provides the unique classical solution to the HJB equation, which is shown to be the value function of our control problem.

\vspace{4mm}

\end{abstract}

\section{Introduction}
The aim of the present paper is to study  optimal control problems for a class of semi-Markov processes using a suitable class of backward stochastic differential equations (BSDEs), driven by the random measure associated to the semi-Markov process itself.

Let us briefly describe our framework.
Our starting point is a semi-Markov pure jump process $X$ on a general
state space $K$.
It is constructed starting from a jump rate function  $\lambda(x,a)$ and a jump measure $A\mapsto\bar{q}(x,a,A)$
on $K$, depending on  $x\in K$ and  $a\ge0 $.
Our approach is to consider a semi-Markov pure jump process as a two dimensional time-homogeneous and strong Markov process $\{(X_s,a_s), \, s \geq 0\}$ with its natural filtration $\mathcal{F}$ and a family of probabilities  $\P^{x,a}$ for $x \in K$, $a \in  [0,\infty)$  such that $\P^{x,a}(X_0=x,a_0=a)=1$.
If the process starts  from $(x,a)$ at
time $t=0$ then the distribution of its first
jump time $T_1$ under $\P^{x,a}$ is described by the formula
\begin{equation}\label{jumptimeintro}
\P^{x,a}(T_1> s)=
    \exp\left(-\int_a^{a+s}\lambda(x,r)\,dr\right),
\end{equation}
and the conditional probability that the process is in $A$ immediately
after a jump at time $T_1=s$ is
$$
  \P^{x,a}(X_{T_1}\in A\,|\,T_1=s ) =
  \bar{q}(x,s,A).
$$
$X_s$ is called the state of the process at time $s$, and $a_s$ is the duration period in this state up to moment $s$:
\begin{eqnarray*}
a_{s}=
\left\{ \begin{array}{ll}
a + s \qquad \qquad \qquad \qquad \qquad \quad \qquad \text{if}\,X_{p} = X_{s} \quad \forall \, 0 \leqslant p \leqslant s,\,\,p,s \in\R,   \\
s - \sup\{\,p : \, 0\leqslant p \leqslant s,\, X_{p} \neq X_{s}\}  \quad  \text{otherwise.}
\end{array} \right.
\end{eqnarray*}
We note that $X$ alone is not a Markov process.
We limit ourselves to the case of a semi-Markov process $X$
such that the survivor function of $T_1$ under $\P^{x,0}$ is absolutely continuous  and admits a hazard rate function $\lambda$
as in \eqref{jumptimeintro}. The holding times of the process
are not necessarily exponentially distributed and can be infinite with
positive probability.
Our main restriction is that the jump rate function $\lambda$ is uniformly bounded,
which implies that the process $X$ is non explosive.

Denoting by $T_n$ the jump times of $X$, we consider
the marked point process $(T_n,X_{T_n})$ and the associated
random measure $p(dt\,dy)= \sum_{n} \delta_{(T_n,X_{T_n})} $
on $(0,\infty)\times K$, where $\delta$
denotes the Dirac measure.
The dual predictable projection $\tilde p$ of $p$ (shortly, the compensator)
has the following explicit expression
$$
\tilde p(ds\,dy)= \lambda(X_{s-},a_{s-})\,\bar{q}(X_{s-},a_{s-},dy)\,ds.
$$

In Section \ref{section semi_markov_HJB} we address an optimal intensity-control problem for the semi-Markov process.
This is formulated in a  classical way by means of a change
of probability measure, see e.g.  \cite{ElK},  \cite{E}, \cite{B}.
We define a class $\cala$ of admissible control processes
$(u_s)_{s \in [0,\,T]}$;
for every fixed $t \in [0,\,T]$ and $(x,a)\in K \times [0,\infty)$,
the cost to be minimized and the corresponding   value function
are
\begin{eqnarray*}
J(t,x,a, u(\cdot)) & = & \sperxaut{
\int_0^{T-t} l(t+ s,X_s,a_s,u_s)\,ds + g(X_{T-t},a_{T-t})},
\\
v(t,x,a)& = & \inf_{u(\cdot)\in\cala }J(t,x,a,u(\cdot)),
\end{eqnarray*}
where $g,l$ are given real functions.
Here $\E_{u,t}^{x,a}$ denotes the  expectation with respect to another probability $\P_{u,t}^{x,a}$,
depending on $t$ and on the control process $u$ and
constructed in such a way that the
compensator  under $\P_{u,t}^{x,a}$ equals
$r(t+s,X_{s-},a_{s-},y,u_s)\, \lambda(X_{s-},a_{s-})\,\bar{q}(X_{s-},a_{s-},dy)\,ds$, for some function $r$
given in advance as another datum of the control
problem.
Since  the process $(X_s,a_s)_{s \geq 0}$ we want to control is time-homogeneous and starts from $(x,a)$ at time $s = 0$, we introduce a temporal translation which allows to define  the cost functional   for all $t \in [0,T]$. For more details see  Remark \ref{rem:controllo 1}.

Our approach to this control problem  consists in
introducing a family of BSDEs parametrized by $(t,x,a)\in [0,T]\times K \times [0,\infty)$:
\begin{equation}\label{intro_BSDE_controllo}
Y^{x,a}_{s,t} + \int_{s}^{T-t}\int_{K}Z^{x,a}_{\sigma,t}(y)\,q(d\sigma\,dy) =
 g(X_{T-t},a_{T-t}) + \int_{s}^{T-t}\,f\Big(t+\sigma,X_{\sigma},a_{\sigma},Z^{x,a}_{\sigma,t}(\cdot)\Big)\,d\sigma, \quad   s\in [0,\,T-t],
\end{equation}
where the generator is given by the Hamiltonian function $f$ defined for every $s \in [0,\,T]$, $(x,a) \in K\times[0,\,+\infty)$, $z \in L^{2}(K,\mathcal{K},\lambda(x,a)\bar{q}(x,a,dy))$, as
\begin{equation}\label{intro_hamilton_function}
f(s,x,a,z(\cdot))=
 \inf_{u \in U}\Big\{\,l(s,x,a,u) + \int_{K}z(y)(r(s,x,a,y,u)-1)\lambda(x,a)\bar{q}(x,a,dy)\,\Big\}.
\end{equation}
Under appropriate assumptions 
we prove that the optimal control
problem has a solution
and that the value function and the optimal
control can be represented by means of the solution to the BSDE \eqref{intro_BSDE_controllo}.

Backward equations driven by random measures
have been studied in many papers, within \cite{TaLi}, \cite{BaBuPa}, \cite{Roy}, \cite{KhMaPhZh}, \cite{Xia},
and more recently \cite{Be}, \cite{Cre-Mat}, \cite{KaTa-Po-Zh_1}, \cite{KaTa-Po-Zh_2}, \cite{CoFu-mpp}, \cite{CoFu-m}.
In many of them  the stochastic
equations are driven by
a Wiener process and  a Poisson process, see, e.g., \cite{TaLi}, \cite{BaBuPa}, \cite{Roy}, \cite{KhMaPhZh}.
A more general
results on BSDEs driven by random measures is given by  \cite{Xia}, but in this case
the generator $f$  depends on the process $Z$ in a specific
way 
and this condition prevents a direct application to optimal control
problems.
In  \cite{Be}, \cite{Cre-Mat}, \cite{KaTa-Po-Zh_1}, \cite{KaTa-Po-Zh_2}, the authors deal with BSDEs with jumps with a random compensator more general than the compensator of a Poisson random measure; here are involved random compensators which are absolutely continuous with respect to a deterministic measure, that can be reduced to a Poisson measure by a Girsanov change of probability.
Finally, BSDEs driven by a random measure  related  to a pure jump process have been recently studied in \cite{CoFu-mpp}, and in \cite{CoFu-m}  the pure jump Markov case is considered.

Our backward equation 
\eqref{intro_BSDE_controllo}
is driven by a random measure associated to a two dimensional Markov process $(X,a)$, and his compensator is a stochastic random measure with a non-dominated intensity as in \cite{CoFu-m}.
Even if the associated process is  not pure jump,
the existence, uniqueness and continuous
dependence on the data for the BSDE \eqref{intro_BSDE_controllo}
can be deduced extending in a straightforward way the
results in \cite{CoFu-m}.

Concerning the optimal control of semi-Markov processes, the case of a finite number of states has been studied in \cite{chito}, \cite{Howard}, \cite{Jewell}, \cite{Osaki}, while the case of arbitrary state space is considered in  \cite{Ross} and  \cite{St1}.
As in \cite{chito} and in \cite{St1}, in our formulation we admit control actions that can depend not only on the state process but also on the length of time the process has remained in that state.
The approach based on BSDEs is classical in the diffusive context and is also present in the literature in the case of BSDEs with jumps, see as instance \cite{LimQuenez}.
However, it seems to us be pursued here for the first time in the case of the semi-Markov processes.
It allows to treat
in a unified way a large class of control problems, where the state space is general and the
running and final cost are not necessarily bounded.
We remark that, comparing with \cite{St1},
the controlled processes we deal with have laws
absolutely continuous with respect to a given, uncontrolled process;
see also a more detailed comment
in Remark \ref{confrontocontrollo} below. Moreover, in \cite{St1} optimal control problems for semi-Markov processes are studied in the case of infinite time horizon.

In Section \ref{section semi_markov_Kolmogorov_equation} we solve a nonlinear variant
of the Kolmogorov equation for the process $(X,a)$, with the BSDEs approach.
The  process $(X,a)$  is time-homogeneous and Markov, but is not a pure jump process. In particular it has the integro-differential infinitesimal generator
\begin{displaymath}
\mathcal{\tilde{L}}\Phi(x,a) := \partial_a \Phi(x,a) + \int_{K}[\Phi(y,0)-\Phi(x,a)]\,\lambda(x,a)\,\bar{q}(x,a,dy), \qquad (x,a)\in K \times [0, \infty).
\end{displaymath}
The additional differential term $\partial_a $ do not allow to study the associated nonlinear Kolmogorov equation proceeding  as in the pure jump Markov processes framework (see \cite{CoFu-m}).
On the other hand,
the two dimensional Markov process $(X_s,a_s)_{s \geqslant 0}$ belongs to the larger class of piecewise-deterministic Markov processes (PDPs)
introduced by M.H.A. Davis in \cite{Da-bo}, and studied
in the optimal control framework by several authors, within
\cite{Da-Fa}, \cite{Ver}, \cite{Dem}, \cite{LenYam}.
Moreover, we deal with a very specific PDP: taking into account the particular structure of semi-Markov processes,
we present a reformulation of the Kolmogorov equation which allows us to consider solutions in a classical sense.
In particular, we notice that the second component of the process $(X_s,\,a_s)_{s \geqslant 0}$ is linear in $s$.
This fact suggests to introduce the formal directional derivative operator
\begin{equation*}
(Dv)(t,x,a):= \lim_{h \downarrow 0}\frac{v(t+h,x,a+h)-v(t,x,a)}{h},
\end{equation*}
and to consider the following nonlinear Kolmogorov equation
\begin{eqnarray}\label{Kolmogorov_diff_rif_intro}
\left\{ \begin{array}{ll}
Dv(t,x,a) +\mathcal{L}v(t,x,a)  + f(t,x,a,v(t,x,a),v(t,\cdot,0)-v(t,x,a))=0 ,\\
\qquad \qquad \qquad \qquad \qquad \qquad \qquad \qquad \qquad \qquad \qquad \qquad t \in [0,T], \, x \in K, \, a \in [0, \infty),  \\
v(T,x,a)= g(x,a),
\end{array} \right.
\end{eqnarray}
where $$\mathcal{L}\Phi(x,a): = \int_{K}[\Phi(y,0)-\Phi(x,a)]\,\lambda(x,a)\,\bar{q}(x,a,dy), \qquad (x,a)\in K \times [0, \infty).$$  Then we look for a  solution $v$ such that the map  $t \mapsto v(t,x,t+c)$ is absolutely continuous on $[0,T]$, for all constants $c \in [-T,\,+ \infty)$. The functions $f,g$ in \eqref{Kolmogorov_diff_rif_intro} are given.
While it is easy to prove well-posedness of
\eqref{Kolmogorov_diff_rif_intro} under boundedness assumptions,
we achieve the purpose of finding a unique solution under much weaker conditions
related to the distribution of the process $(X,a)$: see
Theorem \ref{thm_kolm}.
To this end we need to define a formula of It$\hat{\mbox{o}}$ type, involving the directional derivative operator $D$, for the composition of the process $(X_s,\,a_s)_{s \geqslant 0}$ with functions $v$ smooth enough (see Lemma \ref{Ito formula} below).\\
We construct the solution $v$
by means of a family of BSDEs of the form
\eqref{intro_BSDE_controllo}.
By the  results above there exists a unique solution
$(Y_{s,t}^{x,a}, Z_{s,t}^{x,a})_{s\in [0,\,T-t]}$ and 
 the estimates on the BSDEs
are used to prove well-posedness of
\eqref{Kolmogorov_diff_rif_intro}. As a by-product we also obtain
the  representation formulae
$$
v(t,x,a)=Y_{0,t}^{x,a},
\qquad
    Y_{s,t}^{x,a}=v(t+s,X_s,a_s),
\qquad
Z_{s,t}^{x,a}(y)= v(t+s,y,0)-  v(t+s,X_{s-},a_{s-}),
$$
which are sometimes called, at least in the diffusive case,
non linear Feynman-Kac formulae.\\
Finally we can go back to the original control problem and observe that the associated Hamilton-Jacobi-Bellman equation has
the form \eqref{Kolmogorov_diff_rif_intro}
where $f$
is the Hamiltonian function \eqref{intro_hamilton_function}. By previous results we are  able to identify  the HJB solution $v(t,x,a)$, constructed probabilistically via BSDEs, with the value function.

\section{Notation, preliminaries and basic assumptions}\label{section_notation}

\subsection{Semi-Markov jump processes}\label{subsection_construction_SMP}
We recall the definition of a semi-Markov process, as given, for instance, in \cite{G-S}.
More precisely we will deal with a semi-Markov process with infinite lifetime (i.e. non explosive).

Suppose we are given a measurable space $(K, \mathcal{K})$, a set $\Omega$ and two functions $X: \Omega \times [0,\infty) \rightarrow K$, $a: \Omega \times [0,\infty) \rightarrow [0,\infty)$. For every $t \geq 0$, we denote by $\mathcal{F}_t$ the $\sigma$-algebra $\sigma((X_s,a_s), \, s \in [0,t] )$. We suppose that for every $x \in K$ and $a \in  [0,\infty)$, a probability $\P^{x,a}$ is given on $(\Omega, \mathcal{F}_{[0,\infty)})$ and the following conditions hold.

\begin{enumerate}
\item  $\mathcal{K}$ contains all one-point sets.
$\Delta$ denotes a point not included in $K$.

\item $\P^{x,a}(X_0=x,a_0=a)=1$ for every $x \in K$, $a \in  [0,\infty)$.
\item  For every $s,\,p \geqslant 0$ and $A\in \calk$ the function $(x,\,a)\mapsto \mathbb{P}^{x,a}(X_{s}\in A,\,a_{s}\leqslant p)$  is $\mathcal{K}\otimes \mathcal{B}^+$-measurable.

  \item  For every  $0\le t\le s$, $p \geqslant 0$, and  $A\in \calk$  we have
 $\P^{x,a}(X_s\in A,\,a_{s}\leqslant p\,|\,\calf_{t})=\P^{X_t,a_t}(X_s\in A,\,a_{s}\leqslant p)$, $\P^{x,a}$-a.s.

\item All the trajectories of the process $X$ have right limits when $K$ is given the discrete topology (the one where all subsets are open). This is equivalent to require that
for every $\omega\in\Omega$ and $t\ge 0$ there
exists $\delta>0$ such that $X_s(\omega)=X_t(\omega)$ for $s\in [t,t+\delta]$.

\item All the trajectories of the process $a$ are continuous from the right piecewise linear functions. For every $\omega\in\Omega$, if $[\alpha, \beta)$ is the interval of linearity of $a_{\cdot}(\omega)$ then $a_s(\omega)= a_{\alpha}(\omega) + s -\alpha$ and $X_{\alpha}(\omega)=X_s(\omega)$; if $\beta$ is a discontinuity point of $a_{\cdot}(\omega)$ then $a_{\beta+}(\omega)=0$ and $X_{\beta}(\omega) \neq X_{\beta-}(\omega)$.

\item For every $\omega\in\Omega$ the number of jumps of the trajectory $t \mapsto X_t(\omega)$ is finite on every bounded interval.
\end{enumerate}
$X_s$ is called the \emph{state} of the process at time $s$, $a_s$ is the   \emph{duration period} in this state up to moment $s$. Also we call $X_s$ the \emph{phase} and $a_s$ the \emph{age} or the \emph{time component} of a semi-Markov process. $X$ is
 a non explosive process because of condition 7.
We note, moreover, that the two-dimensional process $(X,a)$ is a strong Markov process with time-homogeneous transition probabilities because of conditions 2,\,3, and 4. It has right-continuous sample paths because of conditions 1, 5 and 6, and it is not a pure jump Markov process, but only a PDP.

The class of semi-Markov processes we  consider in the paper will be described by means of a special form
of joint law $Q$ under $\mathbb{P}^{x,a}$ of the first jump time $T_{1}$, and the corresponding position $X_{T_{1}}$.
To proceed formally, we fix $X_{0} = x \in K$ and define the first jump time
\begin{equation*}
T_{1} = \inf\{ p >0:\, X_{p} \neq x \},
\end{equation*}
with the convention that $T_{1} = +\infty$ if the indicated set is empty.\\
We introduce $S := K  \times [0,\,+\infty) $ an we denote by $\mathcal{S}$  the smallest $\sigma$-algebra containing all sets of $\mathcal{K} \otimes \mathcal{B}([0,\,+\infty))$. (Here and in the following $\calb(\Lambda)$ denotes
the Borel $\sigma$-algebra of a topological space $\Lambda$).
Take an extra point $\Delta \notin K$ and define $X_{\infty}(\omega)= \Delta$ for all $\omega \in \Omega$, so that $X_{T_{1}}: \Omega \rightarrow K \cup \{ \Delta \}$ is well defined.
Then on the extended space $S  \cup \{ (\Delta ,\,\infty)  \}$ we consider the smallest $\sigma$-algebra, denoted by $\mathcal{S}^{\text{enl}}$, containing $ \{ (\Delta ,\,\infty) \}$ and all sets of $\mathcal{K} \otimes \mathcal{B}([0,\,+\infty))$.
Then $(X_{T_{1}},\,T_{1})$ is a random variable with values in $(S\cup  \{ (\Delta ,\,\infty) \},\mathcal{S}^{\text{enl}})$.
Its law under $\mathbb{P}^{x,a}$ will be denoted by $Q(x,a,\cdot)$.

We will assume that $Q$ is constructed   from two given functions denoted by $\lambda$ and $\bar{q}$. More precisely we assume the following.
\begin{hypothesis}\label{hp_dati}
There exist two functions
\begin{equation*}\lambda: S \rightarrow [0,\infty)  \mbox{  and } \bar{q}:S \times \mathcal{K}\rightarrow [0,1] \end{equation*}
such that
\begin{itemize}
\item[(i)] $(x, a) \mapsto \lambda(x,a)$ is $\mathcal{S}$-measurable;
\item[(ii)] $\sup_{(x,a) \in S}\lambda(x,a) \leqslant C \in \R^+$;
\item[(iii)] $(x, a) \mapsto \bar{q}(x,a,A)$ is $\mathcal{S}$-measurable $\forall A \in \mathcal{K}$;
\item[(iv)] $A \mapsto \bar{q}(x,a,A)$ is a probability measure on $\mathcal{K}$ for all $(x,\,a)\in S$.
\end{itemize}
\end{hypothesis}
We define a function $H$ on $ K \times[0,\infty]$ by
\begin{equation}\label{def_H}
H(x, s):= 1- e^{-\int_{0}^{s}\lambda(x,r)dr}.
\end{equation}
Given $\lambda$ and $\bar{q}$, we will require that for the semi-Markov process $X$ we have,
for every $(x,a)\in S$
and for  $A\in \calk$, $0\le c< d\le \infty$,

\begin{eqnarray}\label{jumpkernel}
    Q(x,a,  A \times (c,d) ) & = & \dis \frac{1}{1- H(x,a)}\int_c^d
    \bar{q}(x,s,A) \frac{d }{d\,s} \, H(x,a+s)\,ds\nonumber\\
    & = & \dis \int_c^d
    \bar{q}(x,s,A)\;\lambda (x,a+s)\;
    \exp\left(-\int_a^{a+s}\lambda(x,r)\,dr\right) \,ds,
\end{eqnarray}
where
  $Q$ was described above as the law of $(X_{T_1},T_1)$
under $\P^{x,a}$.

The existence of a semi-Markov process satisfying
\eqref{jumpkernel} is a well known fact, see for instance
\cite{St1} Theorem 2.1, where it is proved
that $X$ is in addition a strong Markov process.
The nonexplosive character of $X$ is made possible by Hypothesis \ref{hp_dati}-(ii).

We note that our data only consist initially in a measurable space
$(K,\calk)$ ($\mathcal{K}$ contains all singleton subsets of $K$), and in two functions $\lambda$, $\bar{q}$ satisfying Hypothesis \ref{hp_dati}.
The semi-Markov process $X$
can be constructed in an arbitrary way provided
\eqref{jumpkernel} holds.

\begin{remark}\label{processosemi-markov}
\begin{enumerate}
\item Note that \eqref{jumpkernel} completely specifies the probability measure  $Q(x,a,\cdot)$ on
    $(S\,\cup\,  \{ (\Delta ,\,\infty) \},\mathcal{S}^{\text{enl}})$: indeed  simple computations show  that, for $s\ge 0$,
    \begin{equation}\label{jumpkerneldue}
     \P^{x,a}(T_1\in (s,\infty])  =1- Q(x,a,K \times (0,s] )
     = \exp\left(-\int_a^{a+s}\lambda(x,r)\,dr\right),
    \end{equation}
    and
    we clearly have
    \begin{equation*}\label{jumpkerneltre}
    \begin{array}{lll}
    \P^{x,a}(T_1=\infty)&=&\dis
    Q(x, a,\{(\Delta,\infty)\})=\exp\left(-\int_a^\infty\lambda(x,r)\,dr\right).
    \end{array}
    \end{equation*}
    Moreover, the kernel $Q$ is well defined, because $H(x,a)< 1$ for all $(x,a) \in S$ by assumption \ref{hp_dati}-(ii).

\item
    The data $\lambda$ and $\bar{q}$ have themselves a probabilistic interpretation.
    In fact  if in  \eqref{jumpkerneldue} we set $a = 0$ we obtain \begin{equation}
    \P^{x,0}(T_1 > s) = \exp\left(-\int_0^s\lambda(x,r)\,dr\right) = 1 - H(x,s).
    \end{equation}
    This means that under $\P^{x,0}$ the law of $T_1$ is described by the distribution function $H$, and
    \begin{displaymath}
    \lambda(x,a)= \frac{\frac{\partial H}{\partial a}(x,a)}{1- H(x,a)}.
    \end{displaymath}
    Then $\lambda(x,a)$ is the jump rate of the process $X$ given that it has been in state $x$ for a time $a$.\\
    Moreover, the probability
    $\bar{q}(x,s,\cdot)$ can be interpreted as the conditional probability that $X_{T_1}$
    is in $A\in\calk$ given that $T_1=s$; more precisely,
    $$
    \P^{x,a}(X_{T_1}\in A, T_1<\infty\,|\,T_1 ) =
    \bar{q}(x,T_1,A)\,1_{T_1<\infty}, \qquad  \P^{x,a}-a.s.
    $$

\item
    In \cite{G-S} the following observation is made: starting from $T_0=t$ define inductively
    $
    T_{n+1}=\inf\{s>T_n\,:\, X_s\neq
    X_{T_{n}}\},
    $
    with the convention that $T_{n+1} =\infty$ if the indicated
    set is empty; then, under the probability
    $\P^{x,a}$, the sequence of the successive states of the semi-Markov $X$ is a Markov chain, as in the case of  Markov processes.
    However, while for the latter the duration period in the state depends only on this state and it is necessarily exponentially distributed, in the case of a semi Markov process the duration period depends also on the state into which the process moves and the distribution of the duration period may be arbitrary.

\item In \cite{G-S} is also
    proved that the sequence $(X_{T_n},T_n)_{n\ge 0}$ is a discrete-time
    Markov process in $(S\cup \{(\Delta,\,\infty)\},$ $\mathcal{S}^{\text{enl}})$ with transition kernel $Q$, provided we extend
    the definition of $Q$ making the state $(\Delta,\,\infty)$ absorbing, i.e. we define
    $$
    Q(\Delta,\,\infty,\, S)=0,\qquad
    Q(\Delta,\,\infty, \,\{(\Delta,\,\infty)\})=1.
    $$
    Note that $(X_{T_n},T_n)_{n\ge 0}$ is time-homogeneous.

    This fact allows for a simple description of the process $X$.
    Suppose one starts with  a  discrete-time Markov process $(\tau_n,\xi_n)_{n\ge 0}$
    in $S$ with transition probability kernel $Q$ and a given starting
    point $(x,a)\in S$ (conceptually, trajectories
    of such a process are easy to simulate). One can then define a process $Y$
     in $K$  setting $Y_t=\sum_{n=0}^N\xi_n 1_{[\tau_n, \tau_{n+1}) }(t)$,
    where $N=\sup\{n\ge 0\,:\, \tau_n\leqslant\infty\}$. Then  $Y$ has the same law
    as the process $X$ under $\P^{x,a}$.

\item We stress that \eqref{def_H} limits ourselves to deal with a class of semi-Markov processes for which the survivor function $T_1$ under $\P^{x,0}$ admits a hazard rate function $\lambda$.

\end{enumerate}
\end{remark}

\subsection{BSDEs driven by a semi-Markov process}
Let be given a measurable space $(K,\mathcal K)$, a transition measure $\bar{q}$ on $K$ and a given positive function $\lambda$, satisfying Hypothesis \ref{hp_dati}. Let  $X$ be the associated semi-Markov process constructed out of them as described in Section \ref{subsection_construction_SMP}.
We fix a deterministic terminal time $T>0$ and a pair $(x,a)\in S$, and we look at all processes under the probability $\mathbb{P}^{x,a}$.
We denote by $\mathcal{F}$ the natural filtration $(\mathcal{F}_t)_{t \in[0, \infty)}$ of $X$.
Conditions 1, 5 and 6  above imply that the filtration $\mathcal{F}$ is right continuous (see \cite{B}, Appendix A2, Theorem T26).
The predictable $\sigma$-algebra (respectively, the progressive $\sigma$-algebra) on $\Omega \times [0, \, \infty)$ is denoted by $\mathcal{P}$ (respectively, by $Prog$). The same symbols  also denote the restriction to $\Omega \times [0, \, T]$.

We define  a sequence $(T_{n})_{n \geqslant 1}$ of random variables with values in $[0, \, \infty]$, setting
\begin{equation}\label{Tn_def}
T_{0}(\omega) = 0,\quad
T_{n+1}(\omega) =\inf\{ s \geqslant T_{n}(\omega): \, X_{s}(\omega) \neq X_{T_{n}}(\omega) \},
\end{equation}
with the convention that $T_{n+1}(\omega)= \infty$ if the indicated set is empty.
Being $X$ a jump process we have $T_{n}(\omega)\leqslant T_{n+1}(\omega)$ if $T_{n+1}(\omega)< \infty$, while the non explosion of $X$ means that $T_{n+1}(\omega)\rightarrow \infty$. We stress the fact that $(T_{n})_{n \geqslant 1}$  coincide by definition with the time jumps of the two dimensional process $(X,a)$.

For $\omega \in \Omega$ we define a random measure on $( [0, \, \infty)\times K ,\, \mathcal{B}[0, \, \infty)\otimes  \mathcal{K})$ setting
\begin{equation}
p(\omega, C) = \sum_{n \geqslant 1}\one_{\{(T_{n}(\omega), \, X_{T_{n}}(\omega)) \in C \}}, \qquad C \in \mathcal{B}[0, \, \infty)\otimes  \mathcal{K}.
\end{equation}
The random measure $\lambda(X_{s-},a_{s-})\,\bar{q}(X_{s-},a_{s-},dy)\,ds$ is called the compensator, or the dual predictable  projection, of $p(ds,dy)$. 
We are interested in
the following family of backward equations driven by the compensated random measure $q(ds\,dy)=p(ds\,dy) -\lambda(X_{s-},a_{s-})\,\bar{q}(X_{s-},a_{s-},dy)\,ds$ and parametrized by $(x,a)$: $\mathbb{P}^{x,a}$-a.s.,
\begin{equation}\label{BSDE}
Y_{s} + \int_{s}^{T}\int_{K}Z_{r}(y)\,q(dr\,dy) =
 g(X_{T},a_{T}) + \int_{s}^{T}f\Big(r,X_{r},a_{r},Y_{r},Z_{r}(\cdot)\Big)\,dr, \qquad s\in [0,\,T].
\end{equation}
We 
consider the following assumptions on the data $f$ and $g$.
\begin{hypothesis}\label{H_1}
\begin{itemize}
\item[\emph{(1)}] The final condition $g : S \rightarrow \R$ is $\mathcal S$-measurable \\
 and $\sperxa{\abs{g(X_{T},a_{T})}^{2}} < \infty$.
\item[\emph{(2)}] The generator $f$ is such that
\begin{itemize}
\item[\emph{(i)}] for every $s \in [0,\,T]$, $(x,a)\in S$, $r \in \R$, $f$ is a mapping \\
$f(s,x,a,r,\cdot) : \mathcal{L}^{2}(K,\mathcal{K},\,\lambda(x,a)\,\bar{q}(x,a,dy))\rightarrow \R$;
\item[\emph{(ii)}] for every bounded and $\mathcal{K}$-measurable $z: K \rightarrow \R$ the mapping
\begin{equation}\label{map}
(s,x,a,r)\mapsto f(s,x,a,r,z(\cdot))
\end{equation}
is $\mathcal{B}([0,\,T]) \otimes \mathcal S \otimes \mathcal{B}(\R)$-measurable;
\item[\emph{(iii)}] there exist $L \geqslant 0$, $L' \geqslant 0 $ such that for every  $s \in [0, \, T]$, $(x,a)\in S$, $r, r' \in \R, \\
 z, z' \in \mathcal{L}^{2}(K,\mathcal{K},\lambda(x,a)\,\bar{q}(x,a,dy))$ we have
\begin{equation}\label{f_inequality}
\abs{f(s,x,a,r,z(\cdot)) - f(s,x,a,r',z'(\cdot))} \leqslant L'\abs{r - r'} + L\left( \int_{K}  \abs{z(y) - z'(y)}^{2} \lambda(x,a)\,\bar{q}(x,a,dy) \right)^{1/2};
\end{equation}
\item[\emph{(iv)}]we have
\begin{equation}\label{finite_sper}
\sperxa{\int_{0}^{T}\abs{f(s,X_{s},a_{s},0,0)}^{2}ds} < \infty.
\end{equation}
\end{itemize}
\end{itemize}
\end{hypothesis}

\begin{remark}
Assumptions (i), (ii), and (iii) imply the following measurability properties of \\
$f(s,X_{s},a_{s},Y_{s},Z_{s}(\cdot))$:
\begin{itemize}
\item if $Z\in \mathcal{L}^2(p)$, then the mapping
\begin{displaymath}
(\omega,s,y) \mapsto f(s,X_{s-}(\omega),a_{s-}(\omega),y,Z_{s}(\omega, \cdot))
\end{displaymath}
is $\mathcal{P}\otimes \mathcal{B}(\R)$-measurable;
\item if, in addition, $Y$ is a $Prog$-measurable process, then
\begin{displaymath}
(\omega,s) \mapsto f(s,X_{s-}(\omega),a_{s-}(\omega),Y_{s}(\omega),Z_{s}(\omega, \cdot))
\end{displaymath}
is $Prog$-measurable.
\end{itemize}
\end{remark}
We introduce the space $\mathbb{M}^{x,a}$ of the processes $(Y,Z)$ on $[0,\,T]$ such that $Y$ is  real-valued and $Prog$-measurable, $Z: \Omega \times K \rightarrow \R$ is $\mathcal{P}\otimes \mathcal{K}$-measurable, and
\begin{displaymath}
||(Y,Z)||^{2}_{\mathbb{M}^{x,a}} := \sperxa{\int_{0}^{T}\abs{Y_{s}}^{2}ds} + \sperxa{\int_{0}^{T}\int_{K}\abs{Z_{s}(y)}^{2}\lambda(X_{s},a_{s})\,\bar{q}(X_{s},a_{s},dy)\,ds} <\infty.
\end{displaymath}
The space $\mathbb{M}^{x,a}$  endowed with this norm is a Banach space, provided we identify pairs of processes whose difference has norm zero.
\begin{theorem}\label{thm: uniqueness_existence_BSDE}
Suppose that Hypothesis \ref{H_1} holds for some $(x,a) \in S$.\\
Then there exists a unique pair $(Y,Z)$ in $\mathbb{M}^{x,a}$ which solves the BSDE \eqref{BSDE}.
Let moreover $(Y',Z')$
be another solution in $\mathbb{M}^{x,a}$ to the BSDE \eqref{BSDE} associated with the driver $f'$  and final datum $g'$. Then
\begin{align}\label{stima-differenza}
&\sup_{s \in [0,\,T]}\sperxa{|Y_s-Y'_{s}|^2} + \sperxa{\int_{0}^{T}|Y_s-Y'_{s}|^2ds}                + \sperxa{\int_{0}^{T}\int_{K}|Z_{s}(y)-Z'_{s}(y)|^2\lambda(X_{s},a_{s})\,\bar{q}(X_{s},a_{s},dy)\,ds} \nonumber\\
& \leqslant C\sperxa{|g(X_{T})-g'(X_{T})|^2 +   \int_{0}^{T}|f(s,X_{s},a_{s},Y'_{s},Z'_{s}(\cdot) )-f'(s,X_{s},a_{s},Y'_{s},Z'_{s}(\cdot))|^2ds},
\end{align}
where $C$ is a constant depending on $T$, $L$, $L'$.
\end{theorem}
\begin{remark}
The construction of a solution to the BSDE \eqref{BSDE} is based on the integral representation theorem of marked point process martingales
(see, e.g., \cite{Da-bo}), and on a fixed-point argument.
Similar results of well-posedness for BSDEs driven by random measures can be found in literature, see, in particular, the theorems given in  \cite{CoFu-m}, Section 3, and  in \cite{Be}.
Notice that  these results can not be a priori straight applied to our framework:
in \cite{Be}  are involved random compensators which are absolutely continuous with respect to a deterministic measure, instead in our case the compensator is a stochastic random measure with a non-dominated intensity;
 \cite{CoFu-m} apply to BSDEs driven by a random measure associated to a pure jump Markov process,  while the two dimensional process $(X,a)$ is Markov but not pure jump.
Nevertheless, under Hypothesis \ref{H_1},
Theorem 3.4 and Proposition 3.5 in \cite{CoFu-m} can be  extended to our framework without additional difficulties.
The proofs turn out to be  very similar to
those of the mentioned results, and we do not report them here to alleviate the presentation.
\end{remark}

\section{Optimal control}\label{section semi_markov_HJB}
\subsection{Formulation of the problem}
In this section we consider again a measurable space $(K,\mathcal{K})$, a transition measure $\bar{q}$ and a function $\lambda$ satisfying Hypothesis \ref{hp_dati}. 
The data specifying the optimal control problem we will address to are an action (or decision) space $U$, a running cost function $l$, a terminal cost function $g$, a (deterministic, finite) time horizon $T>0$ and another function $r$ specifying the effect of the control process.
We define  an admissible control process, or simply a control, as a predictable process $(u_{s})_{s \in [0, \, T]}$ with values in $U$. The set of admissible control processes is denoted by $\mathcal{A}$.
We will make the following assumptions:
\begin{hypothesis} \label{hyp_control}
\begin{itemize}
\item[\emph{(1)}] $(U,\mathcal U)$ is a measurable space.
\item[\emph{(2)}] The function $r: [0, \, T]\times S \times K \times U \rightarrow \R $ is
$\mathcal{B}([0,\,T]) \otimes \mathcal S \otimes \mathcal K \otimes \mathcal U$-measurable and there exists a constant $C_{r}>1$ such that,
\begin{equation}\label{hyp: r_bound}
0\leqslant r(t,x,a,y,u)\leqslant C_{r}, \qquad t\in[0, \, T],\,(x,a)\in S,\, y \in K,\, u \in U.
\end{equation}
\item[\emph{(3)}] The function $g: S \rightarrow \R $ is
$ \mathcal S$-measurable, and for all fixed $t \in [0,\,T]$,
\begin{equation}\label{hyp: g_boundedness}
\sperxa{\abs{g(X_{T-t},a_{T-t})}^{2}}<\infty,\qquad \forall (x,a)\in S.
\end{equation}
\item[\emph{(4)}]  The function $l: [0, \, T]\times S\times U \rightarrow \R $ is
$\mathcal{B}([0\,\,T]) \otimes \mathcal S \otimes \mathcal U$-measurable and there exists $\alpha >1$ such that, for every fixed $t \in [0,\,T]$, for every $(x,a) \in S$ and $u(\cdot) \in \mathcal A$,
\begin{eqnarray}\label{hyp: l_inf_bound}
\begin{array}{ll}
\inf_{u \in U}l(t,x,a,u) > \infty;\\
\sperxa{\int_{0}^{T-t}\abs{\inf_{u \in U}l(t+s,X_{s},a_{s},u)}^{2}\,ds}<\infty,\\
\sperxa{\int_{0}^{T-t}\abs{l(t+s,X_{s},a_{s},u_s)}\,ds}^{\alpha}<\infty.
\end{array}
\end{eqnarray}
\end{itemize}
\end{hypothesis}
To any $(t,x,a) \in [0,\,T] \times S$ and any control $u(\cdot)\in \mathcal{A}$ we associate a probability measure $\mathbb{P}^{x,a}_{u,t}$ by a change of measure of Girsanov type, as we now describe.
Recalling the definition of the jump times $T_{n}$ in \eqref{Tn_def}, we define, for every fixed $t \in [0,\,T]$,
$$
L^{t}_{s} = \exp\left(\int_{0}^{s}\!\!\int_{K}(1-r(t+\sigma,X_{\sigma},a_{\sigma},y,u_{\sigma}) )\,\lambda(X_{\sigma},a_{\sigma})\,\bar{q}(X_{\sigma},a_{\sigma},dy)\,d\sigma \!\right)\!\!\!\prod_{n\geqslant 1: T_{n} \leqslant s}\!\!\!r(t+T_{n}, X_{T_{n}-},a_{T_{n}-},X_{T_{n}},u_{T_{n}}),
$$
for all $s \in [0,\,T-t]$,
with the convention that the last product equals $1$ if there are no indices $n \geqslant 1$ satisfying $T_{n}\leqslant s$.
As a consequence of the boundedness assumption on $\bar{q}$ and $\lambda$ it can be proved, using for instance Lemma 4.2 in \cite{CoFu-mpp}, or \cite{B} Chapter VIII  Theorem T11, that for every fixed $t \in [0,\,T]$ and for every $\gamma >1$ we have
\begin{equation}\label{L_martingale}
\sperxa{\abs{L^t_{T-t}}^{\gamma}}<\infty,\qquad \sperxa{L^t_{T-t}}=1,
\end{equation}
and therefore the process $L^t$ is a martingale (relative to $\mathbb{P}^{x,a}$ and $\mathcal{F}$). Defining a probability $\mathbb{P}_{u,t}^{x,a}(d\omega)= L^t_{T-t}(\omega)\,\mathbb{P}^{x,a}(d\omega)$, we introduce the cost functional corresponding to $u(\cdot)\in \mathcal{A}$ as
\begin{equation}\label{cost_functional}
J(t,x,a,u(\cdot))= \sperxaut{ \int_{0}^{T-t}\,l(t+s,X_{s},a_{s},u_{s})\,ds + g(X_{T-t},a_{T-t}) },
\end{equation}
where $\mathbb{E}_{u,t}^{x,a}$ denotes the expectation under $\mathbb{P}_{u,t}^{x,a}$.
Taking into account \eqref{hyp: g_boundedness}, \eqref{hyp: l_inf_bound} and \eqref{L_martingale}, and using H\"{o}lder inequality it is easily seen that the cost is finite for every admissible control. The control problem starting at $(x,a)$ at time $s=0$ with terminal time $s=T-t$ consists in minimizing $J(t,x,a,\cdot)$ over $\mathcal{A}$.

We finally introduce the value function
\begin{displaymath}
v(t,x,a) = \inf_{u(\cdot)\in \mathcal{A}}J(t,x,a,u(\cdot)),\qquad t \in [0,\,T],\,\, (x,a)\in S.
\end{displaymath}
The previous formulation of the optimal control problem by means of change of probability measure is classical (see e.g. \cite{ElK}, \cite{E}, \cite{B}). Some comments may be useful at this point.
\begin{remark}\label{rem:controllo 1}
\begin{itemize}
\item[1.] The particular form of cost functional \eqref{cost_functional} is due to the fact that the time-homogeneous Markov process ${(X_{s},a_{s})}_{s \geqslant 0}$ satisfies
    \begin{equation*}
    \P^{x,a}(X_{0}=x,\,a_{0}=a) = 1;
    \end{equation*}
    the introduction of the temporal translation in the first component allows us to define $J(t,x,a, u(\cdot))$ for all $t \in [0,\,T]$.
\item [2.] We recall (see e.g. \cite{B}, Appendix A2, Theorem T34) that a process $u$ is $\mathcal{F}$-predictable if and only if it admits the representation
    \begin{displaymath}
    u_{s}(\omega) = \sum_{n \geqslant 0}u_{s}^{(n)}(\omega)\,\one_{(T_{n}(\omega),T_{n+1}(\omega)]}(s)
    \end{displaymath}
    where for each $(\omega,s)\mapsto u_{s}^{(n)}(\omega)$ is $\mathcal{F}_{[0,\,T_{n}]}\otimes \mathcal{B}(\R^+)$-measurable, with $\mathcal{F}_{[0,\,T_{n}]}= \sigma(T_{i},X_{T_{i}},\, 0\leqslant i \leqslant n)$ (see e.g. \cite{B}, Appendix A2, Theorem T30).
    Thus the fact that controls are predictable processes admits the following interpretation: at each time $T_{n}$ (i.e. immediately after a jump) the controller, having observed the random variables $T_{i},\, X_{T_{i}},\, (0\leqslant i \leqslant n)$, chooses his current action, and updates her/his decisions only at time $T_{n+1}$.
\item[3.] It can be proved (see \cite{J} Theorem 4.5) that the compensator of $p(ds \,dy)$ under $\mathbb{P}_{u,t}^{x,a}$ is
    \begin{displaymath}
    r(t+s,X_{s-},a_{s-},y,u_{s})\,\lambda(X_{s-},a_{s-})\,\bar{q}(X_{s-},a_{s-},dy)\,ds,
    \end{displaymath}
    whereas the compensator of  $p(ds\, dy)$ under $\mathbb{P}^{x,a}$ was $\lambda(X_{s-},a_{s-})\,\bar{q}(X_{s-},a_{s-},dy)\,ds$. This explains that the choice of a given control $u(\cdot)$ affects the stochastic system multiplying its compensator by $r(t+s,x,a,y,u_{s})$.
\item[4.]We call control law an arbitrary measurable function $\underbar{u}: [0,\,T]\times S\rightarrow U$. Given a control law one can define an admissible control $u$ setting $u_{s}= \underbar{u}(s,X_{s-},a_{s-})$.\\ Controls of this form are called feedback controls. For a feedback control the compensator of $p(ds \,dy)$ is
     $r(t+s,X_{s-},a_{s-},y,\underbar{u}(s,X_{s-},a_{s-}))\,\lambda(X_{s-},a_{s-})\,\bar{q}(X_{s-},a_{s-},dy)\,ds$ under $\mathbb{P}_{u,t}^{x,a}$. Thus, the process $(X,a)$ under the optimal probability is a two-dimensional Markov  process corresponding to the transition measure
    \begin{displaymath}
    r(t+s,x,a,y,\underbar{u}(s,x,a))\,\lambda(x,a)\,\bar{q}(x,a,dy)
    \end{displaymath}
    instead of $\lambda(x,a)\,\bar{q}(x,a,dy)$.
    However, even if the optimal control is in the feedback form, the optimal process is not, in general, time-homogeneous since the control law may depend on time. In this case, according to the definition given in Section \ref{section_notation}, the process $X$ under the optimal probability is not a semi-Markov process.
\end{itemize}
\end{remark}

\begin{remark}\label{confrontocontrollo}
Our formulation of the optimal control should be compared with another approach (see e.g. \cite{St1}).
In \cite{St1} is given a family of jump measures on K $\{\bar{q}(x, b ,\cdot),\,b \in B\}$ with $B$ some index set endowed with a topology.
In the so called \emph{strong formulation} a control $u$ is an ordered pair of functions $(\lambda',\beta)$ with $\lambda': S \rightarrow \R^+$, $\beta: S \rightarrow B$ such that
\begin{align*}
\begin{array}{lll}
\lambda'\,\,\text{ and}\,\, \beta\,\, \text{are}\,\, \mathcal{S}-\text{measurable}; \\
\forall x \in K,\,	\exists\, t(x)>0:\,\,\int_{0}^{t(x)}\,\lambda'(x,r)\, dr < \infty;\\
\bar{q}(\cdot, \beta ,A)\,\, \text{is}\,\, \mathcal{B}^+\text{-measurable}\,\,  \forall A \in \mathcal{K}.
\end{array}
\end{align*}
If $\mathcal{A}$ is the class of controls which satisfies the above conditions, then a control $u = (\lambda', \beta)\in \mathcal{A}$ determines a controlled process $X^{u}$ in the following manner.
Let
\begin{displaymath}
H^{u}(x, s):= 1- e^{-\int_{0}^{s}\,\lambda'(x,r)\,dr}, \quad \forall (x,s)\in S,
\end{displaymath}
and suppose that $(X^{u}_{0},a_{0}^{u})=(x,a)$.
Then at time $0$, the process starts in state $x$ and remains there a random time $S_{1}>0$, such that
\begin{equation}\label{tempo_soggiorno_controllato}
\probxa{S_{1} \leqslant s} = \frac{H^{u}(x,a+s)-H^{u}(x,a)}{1- H^{u}(x,a)}.
\end{equation}
At time $ S_{1}$ the process transitions to the state $X^{u}_{ S_{1}}$, where
\begin{displaymath}
\probxa{X^{u}_{S_{1}} \in A| S_{1}} = \bar{q}(x,\beta(x,S_{1}),A).
\end{displaymath}
The process stays in state $X^{u}_{ S_{1}}$ for a random time $S_{2}>0$ such that
\begin{displaymath}
\probxa{S_{2} \leqslant s| S_{1}, \,X^{u}_{S_{1}}} = H^{u}(X^{u}_{ S_{1}},s)
\end{displaymath}
and then at time $ S_{1}+S_{2}$ transitions to $X^{u}_{ S_{1}+S_{2}}$, where
\begin{displaymath}
\probxa{X^{u}_{ S_{1}+S_{2}} \in A| S_{1},\, X^{u}_{S_{1}},\,S_{2}} = \bar{q}(X^{u}_{ S_{1}},\beta(X^{u}_{S_{1}},S_{2}),A).
\end{displaymath}
We remark that the process $X^u$ constructed in this way turns out to be semi-Markov.

We also mention that the class of control problems specified by the initial data  $\lambda'$ and $\beta$ is in general larger that the one we address in this paper. This can be seen noticing that in our framework all the controlled processes have laws which are absolutely continuous with respect to a single uncontrolled process (the one corresponding to $r\equiv 1$) whereas this might not be the case for the rate measures $\lambda'(x,a)\,\bar{q}(x,\beta(x,a),A)$ when $u= (\lambda',\,\beta)$ ranges in the set of all possible control laws.
\end{remark}

\subsection{BSDEs and the synthesis of the optimal control}
We next proceed to solve the optimal control problem formulated above.
A basic role is played by the BSDE: for every fixed $t \in [0,\,T]$, $\mathbb{P}^{x,a}$-a.s.
\begin{equation}\label{BSDE_controllo}
Y^{x,a}_{s,t} + \int_{s}^{T-t}\int_{K}Z^{x,a}_{\sigma,t}(y)q(d\sigma\,dy) =
 g(X_{T-t},a_{T-t}) + \int_{s}^{T-t}f\Big(t+\sigma,X_{\sigma},a_{\sigma},Z^{x,a}_{\sigma,t}(\cdot)\Big)d\sigma, \quad  \forall s\in [0,\,T-t],
\end{equation}
with terminal condition  given by the terminal cost $g$ and generator given by the Hamiltonian function $f$ defined for every $s \in [0,\,T],\, (x,a) \in S,\, z \in L^{2}(K,\mathcal{K},\,\lambda(x,a)\,\bar{q}(x,a,dy))$, as
\begin{equation}\label{hamilton_function}
f(s,x,a,z(\cdot))=
 \inf_{u \in U}\Big\{\,l(s,x,a,u) + \int_{K}z(y)(r(s,x,a,y,u)-1)\lambda(x,a)\bar{q}(x,a,dy)\,\Big\}.
\end{equation}
In \eqref{BSDE_controllo} the superscript $(x,a)$ denotes the starting point at time $s=0$ of the process $(X_s,\,a_s)_{s \geqslant 0}$, while the dependence of $Y$ and $Z$ on the  parameter $t$ is related to the temporal horizon of the considered optimal control problem.
For every $t \in [0\,\,T]$, we look for a process $Y^{x,a}_{s,t}(\omega)$  adapted and càdlàg and a process $Z^{x,a}_{s,t}(\omega,y)$  $\mathcal{P}\otimes \mathcal{K}$-measurable  satisfying the   integrability conditions
$$
\sperxa{\int_{0}^{T-t}\abs{Y^{x,a}_{s,t}}^{2}ds}<\infty,\qquad
\sperxa{\int_{0}^{T-t}\int_{K}\abs{Z^{x,a}_{s,t}(y)}^{2}\lambda(X_{s},a_{s})\,\bar{q}(X_{s},a_{s},dy)\,ds}<\infty.
$$
One can verify that, under Hypothesis \ref{hyp_control} on the optimal control problem,  all the assumptions of Hypothesis \ref{H_1} hold true for the generator $f$ and the terminal condition $g$ in the BSDE \eqref{BSDE_controllo}.
The only non trivial verification is the Lipschitz condition \eqref{f_inequality}, which follows from the boundedness assumption \eqref{hyp: r_bound}. Indeed, for every $s\in [0,\,T]$, $(x,a)\in S$, $z,\,z'\in L^{2}(K, \mathcal{K},\lambda(x,a)\,\bar{q}(x,a,dy))$,
\begin{align*}
&\int_{K}z(y)(r(s,x,a,y,u))-1)\,\lambda(x,a) \, \bar{q}(x,a,dy)\\
& \leqslant \int_{K}\abs{z(y)-z'(y)}\,(r(s,x,a,y,u)-1)\,\lambda(x,a) \, \bar{q}(x,a,dy) + \int_{K}z'(y)(r(s,x,a,y,u)-1)\,\lambda(x,a) \, \bar{q}(x,a,dy)\\
&  \leqslant (C_{r} + 1)\,(\lambda(x,a) \, \bar{q}(x,a,K))^{1/2}\,\cdot\left(\int_{K}\abs{z(y)-z'(y)}^2\,\lambda(x,a) \, \bar{q}(x,a,dy)\right)^{1/2}\\
& + \int_{K}z'(y)(r(s,x,a,y,u)-1)\,\lambda(x,a) \, \bar{q}(x,a,dy),
\end{align*}
so that, adding  $l(s,x,a,u)$ on both sides and taking the infimum over $u \in U$, it follows that
\begin{equation}
f(s,x,a,z)\leqslant L\left(\int_{K}\abs{z(y)-z'(y)}^2\lambda(x,a) \, \bar{q}(x,a,dy)\right)^{1/2}
 + f(s,x,a,z'),
\end{equation}
where $L:= (C_{r} + 1)\sup_{(x,a) \in S}\,(\lambda(x,a)\,\bar{q}(x,a,K))^{1/2}$;
exchanging $z$ and $z'$ roles we obtain  \eqref{f_inequality}.

Then by Theorem \ref{thm: uniqueness_existence_BSDE}, for every fixed $t \in [0,\,T]$, for every $(x,a)\in S$, there exists a unique solution of \eqref{BSDE_controllo} $(Y^{x,a}_{s,t},Z^{x,a}_{s,t})_{s \in [0,\,T-t]}$, and $Y_{0,t}^{x,a}$ is deterministic. Moreover, we have the following result:
\begin{proposition}
Assume that Hypotheses  \ref{hyp_control} hold. Then, for every $t \in [0,\,T]$, $(x,a) \in S$, and for every $u(\cdot) \in \mathcal{A}$,
\begin{displaymath}
Y_{0,t}^{x,a} \leqslant J(t,x,a,u(\cdot)).
\end{displaymath}
\end{proposition}

\proof
We consider the BSDE \eqref{BSDE_controllo} at time $s=0$ and we apply the expected value $\mathbb{E}_{u,t}^{x,a}$ associated to the controlled probability $\mathbb{P}_{u,t}^{x,a}$.
Since the $\mathbb{P}_{u,t}^{x,a}$-compensator of $p(ds dy)$ is\\
$r(t+s,X_{s-},a_{s-},y,u_{s})\,\lambda(X_{s-},a_{s-})\,\bar{q}(X_{s-},a_{s-},dy)\,ds$, we have that
\begin{align*}
\sperxaut{\int_{0}^{T-t}\int_{K}Z_{s,t}^{x,a}(y)\,q(ds dy)} &= \sperxaut{\int_{0}^{T-t}\int_{K}Z_{s,t}^{x,a}(y)\,p(ds dy)}  \\
&\quad -\sperxaut{\int_{0}^{T-t}\int_{K}Z_{s,t}^{x,a}(y)\,\lambda(X_{s},a_{s})\,\bar{q}(X_{s},a_{s},dy)\,ds}\\
& = \sperxaut{\int_{0}^{T-t}\!\!\int_{K}Z_{s,t}^{x,a}(y)\,[ r(t+s,X_{s},a_{s},y,u_{s})-1 ]\,
\lambda(X_{s},a_{s})\,\bar{q}(X_{s},a_{s},dy)\,ds}\!.
\end{align*}
Then
\begin{align*}
Y_{0,t}^{x,a} &= \sperxaut{g(X_{T-t},a_{T-t})} + \sperxaut{\int_{0}^{T-t}f(t+s,X_{s},a_{s},Z_{s,t}^{x,a}(\cdot))\,ds}\\
& - \sperxaut{\int_{0}^{T-t}\int_{K}Z_{s,t}^{x,a}(y)\,[r(t+s,X_{s},a_{s},y,u_{s})-1]\,
\lambda(X_{s},a_{s})\,\bar{q}(X_{s},a_{s},dy)\,ds}.
\end{align*}
Adding and subtracting $\sperxaut{\int_{0}^{T-t}l(t+s, X_{s},a_{s},u_{s})\,ds}$ on the right side we obtain the following relation:
\begin{align}\label{BSDE_fundam_rel}
&Y_{0,t}^{x,a} = J(t,x,a, u(\cdot))+ \sperxaut{\int_{0}^{T-t}
\left[f(t+s,X_{s},a_{s},Z_{s,t}^{x,a}(\cdot))
-l(t+s, X_{s},a_{s},u_{s})\right]\,ds} \nonumber\\
&\qquad  -\sperxaut{\int_{0}^{T-t}\int_{K} Z_{s,t}^{x,a}(\cdot)\,[r(t+s,X_{s},a_{s},y,u_{s})-1]\,
\lambda(X_{s},a_{s})\,\bar{q}(X_{s},a_{s},dy)\,ds}.
\end{align}
By the definition of the Hamiltonian function $f$, the two last terms are non positive, and it follows that
\[
Y_{0,t}^{x,a} \leqslant J(t,x,a, u(\cdot)),\qquad \forall u(\cdot)\in \mathcal{A}.
\]
\endproof
\noindent We define the following, possibly empty, set:
\begin{align}\label{Gamma}
&\Gamma(s,x,a,z(\cdot))= \{ \,u \in U: f(s,x,a,z(\cdot))=l(s,x,a,u) + \int_{K}z(y)\,(r(s,x,a,y,u)-1)\,\lambda(x,a)\,\bar{q}(x,a,dy);\nonumber\\
& \qquad  \qquad \qquad \qquad \qquad \qquad \qquad  \,s \in [0,\,T],\,(x,a)\in S,\,z \in L^{2}(K,\mathcal{K},\lambda(x,a)\,\bar{q}(x,a,dy))\,\}.
\end{align}
In order to prove the existence of an optimal control we need to require that the infimum in the definition of $f$ is achieved. Namely we assume that
\begin{hypothesis}\label{hyp_assumed_min}
The sets $\Gamma$ introduced in \eqref{Gamma} are non empty; moreover, for every fixed $t \in [0,\,T]$ and  $(x,a)\in S$, one can find an $\mathcal{F}$-predictable process $u^{\ast\,t,x,a}(\cdot)$ with values in $U$ satisfying
\begin{equation}\label{u_in_Gamma}
u^{\ast\,t,x,a}_{s} \in \Gamma(t+s, X_{s-}, a_{s-},Z^{x,a}_{s,t}(\cdot)), \qquad \mathbb{P}^{x,a}-\text{a.s}. \,\,\forall s \in [0,\,T-t].
\end{equation}
\end{hypothesis}
\begin{theorem}\label{control_sol_BSDEs}
Under Hypothesis \ref{hyp_control} and \ref{hyp_assumed_min} for every fixed $t \in [0,\,T]$ and  $(x,a) \in S$,  $u^{\ast\,t,x,a}(\cdot) \in \mathcal{A}$ is an optimal control for the control problem starting from $(x,a)$ at time $s=0$ with terminal value $s=T-t$. Moreover,  $Y_{0,t}^{x,a}$ coincides with the value function, i.e.  $Y_{0,t}^{x,a}=J(t,x,a,u^{*\,t,x,a}(\cdot))$.
\end{theorem}
\proof
It follows immediately from the relation \eqref{BSDE_fundam_rel} and from the definition of the Hamiltonian function $f$.
\endproof
We recall that general conditions can be formulated for the existence of a process $u^{\ast\,t,x,a}(\cdot)$ satisfying \eqref{u_in_Gamma}, hence of an optimal control;  this is done  by means of an appropriate selection theorem, see e.g. Proposition 5.9 in \cite{CoFu-m}.

We end this section with an example where the BSDE \eqref{BSDE_controllo} can be explicitly solved and a closed form solution of an optimal control problem can be found.
\begin{example}
We consider a fixed time interval $[0,\,T]$ and a state space consisting of three states: $K = \{x_1,x_2,x_3,x_4\}$.
We introduce $(T_n,\xi_n)_{n \geqslant 0}$ setting $(T_0,\xi_0)=(0,x_1)$, $(T_n,\xi_n)=( + \infty, x_1)$ if $n \geqslant 3$ and on $(T_1,\xi_1)$ and $(T_2,\xi_2)$ we make the following assumptions: $\xi_1$ takes values $x_2$ with probability $1$, $\xi_2$ takes values $x_3,x_4$ with probability $1/2$. This means that the
system  starts at time zero in a given state $x_1$, jumps into state $x_2$ with probability $1$  at the random time $T_1$ and into state $x_3$ or $x_4$ with equal probability  at the random time $T_2$. It has no jumps after.
We take $U=[0,\,2]$ and define the function $r$ specifying the effects of the control process as $r(x_1,u)=r(x_2,u)=1$, $r(x_3,u)=u$, $r(x_4,u)=2-u$, $u \in U$.
Moreover, the final cost $g$ assumes the value  $1$ in $(x,a)=(x_4,T-T_2)$
and zero otherwise, and the running cost is defined as $l(s,x,a,u)= \frac{\alpha\,u}{2} \,\lambda(x,a)$, where $\alpha >0$ is a fixed parameter.
The BSDE we want to solve takes the form:
\begin{equation}\label{example}
Y_s + \int_{s}^{T}\int_K Z_{\sigma}(y) \,p(d\sigma\,dy) = g(X_T,\,a_T)+\int_{s}^{T}\inf_{u \in [0,\,2]}\left \{ \frac{\alpha \,u }{2} + \int_K Z_{\sigma}(y)\,r(y,u)\,\bar{q}(X_{\sigma},a_{\sigma},dy)\right\} \lambda(X_{\sigma},a_{\sigma})d\sigma
\end{equation}
that can be written as
\begin{align*}
Y_s + 
\sum_{n \geqslant 1} Z_{T_n}(X_{T_n}) \,\one_{\{s < T_n \leqslant T\}}  & =  g(X_T,\,a_T)+\int_{s}^{T}\inf_{u \in [0,\,2]}\left \{ \frac{\alpha \,u }{2}+ Z_{\sigma}(x_2) \right\}\lambda(x_1,a+\sigma) \one_{\{0 \leqslant \sigma  <T_1 \wedge T\}} \,d\sigma\\
  & +  \int_{s}^{T}\inf_{u \in [0,\,2]}\left \{ \frac{\alpha \,u }{2}+ Z_{\sigma}(x_3)\frac{u}{2} + Z_{\sigma}(x_4)(1-\frac{u}{2}) \right\}\lambda(x_2,\sigma-T_1)  \one_{\{T_1 \leqslant\sigma  <T_2 \wedge T\}}\,d\sigma.
\end{align*}
It is known by \cite{CFJ} that BSDEs of this type   admit the following explicit solution $(Y_s,Z_s(\cdot))_{s\in[0,\,T]}$:
\begin{eqnarray*}
Y_s &=&  y^0(s)\one_{\{ s  <T_1\}}+y^1(s,T_1,\xi_1)\,\one_{\{T_1 \leqslant s  <T_2\}}+y^2(s,T_2,\xi_2,T_1,\xi_1)\,\one_{\{ T_2 \leqslant s\}}\\
Z_s(y) &= &  z^0(s,y)\,\one_{\{ s  \leqslant T_1\}}+z^1(s,y,T_1,\xi_1)\,\one_{\{T_1 < s  \leqslant T_2\}},\quad y \in K.
\end{eqnarray*}
To deduce $y^0$ and  $y^1$ we reduce the BSDE to a system of two ordinary differential equation. To this end, it suffices to consider the following cases:
\begin{itemize}
\item{$\omega \in \Omega \,\,\text{such that}\,\, T < T_1(\omega) < T_2(\omega)$:} \eqref{example} reduces to
\begin{align}\label{y0}
y^0(s) & = \int_{s}^{T}\inf_{u \in [0,\,2]}\left \{ \frac{\alpha \,u }{2}+ z^0(\sigma,x_2) \right\}\lambda(x_1,a+\sigma) \,d\sigma
 =  \int_{s}^{T} z^0(\sigma,x_2)\, \lambda(x_1,a+\sigma) \,d\sigma\nonumber\\
& =  \int_{s}^{T}(y^1(\sigma,\sigma,x_2)-y^0(\sigma))\,\lambda(x_1,a+\sigma) \,d\sigma;
\end{align}
\item{$\omega \in \Omega \,\,\text{such that}\,\, T_1(\omega) < T < T_2(\omega)$, \,$s > T_1$:} \eqref{example} reduces to
\begin{eqnarray}\label{y1}
y^1(s, T_1,\xi_1) & = & \int_{s}^{T}\inf_{u \in [0,\,2]}\left \{ \frac{\alpha \,u }{2}+ z^1(\sigma,x_3,T_1,\xi_1)\frac{u}{2} + z^1(\sigma,x_4,T_1,\xi_1)(1-\frac{u}{2}) \right\}\lambda(\xi_1,\sigma-T_1) \,d\sigma
\nonumber\\
& = & \int_{s}^{T}[z^1(\sigma,x_4,T_1,\xi_1)\wedge (\alpha + z^1(\sigma,x_3,T_1,\xi_1))]\,\lambda(\xi_1,\sigma-T_1) \,d\sigma\nonumber\\
& = & \int_{s}^{T}[(1\wedge \alpha)-y^1(\sigma,T_1,\xi_1)]\,\lambda(\xi_1,\sigma-T_1)\, d\sigma.
\end{eqnarray}
\end{itemize}
Solving \eqref{y0} and \eqref{y1} we obtain
\begin{align*}
&y^0(s) = (1 \wedge \alpha)\left( 1- e^{-\int_s^T\lambda(x_1,a+\sigma)\,d\sigma} \right) - (1 \wedge \alpha)\,e^{- \int_s^T\lambda(x_1,a+\sigma)\,d\sigma}\int_s^T \lambda(x_1,a+\sigma)\, e^{\int_{\sigma}^T\lambda(x_1,a+z)\,dz } e^{-\int_{\sigma}^{T}\lambda(x_2,z - \sigma)\,dz}\,d\sigma\}, \\
&y^1(s,T_1, \xi_1) = (1 \wedge \alpha)\left(1 - e^{-\int_{s}^{T}\lambda(\xi_1,\sigma - T_1)\,d\sigma} \right);
\end{align*}
moreover,
\begin{align*}
& y^2(s,T_2, \xi_2, T_1,\xi_1)=\one_{\{\xi_2 = x_4\}},\\
&z^0(s,x_1)=z^0(s,x_3)=z^0(s,x_4)=0,\,\,\, \quad \quad  z^0(s,x_2)=y^1(s,s,x_2)-y^0(s),\\
& z^1(s,x_1,T_1,\xi_1)= z^1(s,x_2,T_1,\xi_1)=0, \quad \quad z^1(s,x_3,T_1,\xi_1) = (1 \wedge \alpha)\left(e^{-\int_{s}^{T}\lambda(\xi_1,\sigma-T_1)\,d\sigma}-1 \right),\\
& z^1(s,x_4,T_1,\xi_1)=1+z^1(s,x_3,T_1,\xi_1),
\end{align*}
where $z^0$ and $z^1$ are obtained respectively from $y^2$, $y^1$ and  $y^1$, $y^0$ by subtraction.\\
\noindent The optimal cost is then given by $Y_0= y^0(0)$. The optimal control is obtained during the computation of the Hamiltonian function:
it is the process
$u_s =
2\,\one _{(T_1, T_2]}(s)$ if $\alpha \leqslant 1$, and the process $u_s=0$
if $\alpha \geqslant 1$  (both are optimal if $\alpha = 1$).
\end{example}

\section{Nonlinear variant of Kolmogorov equation}\label{section semi_markov_Kolmogorov_equation}
Throughout this section we still assume that a semi-Markov process $X$ is given. It is constructed as in Section \ref{subsection_construction_SMP} by the rate function $\lambda$ and the measure $\bar{q}$ on $K$, and $(X,a)$ is the associated time-homogeneous Markov process.
We assume that $\lambda$ and $\bar{q}$ satisfy Hypothesis \ref{hp_dati}.

It is our purpose to present here some nonlinear variants of the classical backward Kolmogorov equation associated to the Markov process $(X,a)$ and 
to show that their solution can be represented probabilistically by means of an appropriate BSDE of the type considered above.

We will suppose that two functions $f$ and $g$ are given, satisfying Hypothesis \ref{H_1}, and that moreover $g$ verifies, for every fixed $t \in [0,\,T]$,
\begin{equation}\label{g_integrability}
\sperxa{\abs{g(X_{T-t},a_{T-t})}^2}< \infty.
\end{equation}
We define the operator
\begin{equation}\label{L_operator}
\mathcal{L}\psi(x,a):= \int_{K}[\psi(y,0)-\psi(x,a)] \,\lambda(x,a)\,\bar{q}(x,a,dy), \qquad (x,a) \in S,
\end{equation}
for every measurable function $\psi : S\rightarrow \R$ for which the integral is well defined.\\
The equation
\begin{align}\label{Kolmogorov_int}
&v(t,x,a)= g(x,a+T-t) + \int_{t}^{T}\mathcal{L}v(s,x,a+s-t)\,ds \\
& + \int_{t}^{T} f(s,x,a+s-t,v(s,x,a+s-t),v(s,\cdot,0)-v(s,x,a+s-t))\,ds, \quad t \in [0,\,T],\,\,(x,a)\in S, \nonumber
\end{align}
with unknown function $v:[0,\,T]\times S \rightarrow \R$ will be called the nonlinear Kolmogorov equation.

Equivalently, one requires that for every $x \in K$ and  for all constant $c \in [-T,\,+ \infty)$,
\begin{eqnarray}
\begin{array}{ll}\label{absolute_continuity}
t \mapsto v(t,x,t+c) \,\,\text{is absolutely continuous on $[0,T]$,}
\end{array}
\end{eqnarray}
and
\begin{align}\label{Kolmogorov_diff}
\left\{ \begin{array}{ll}
Dv(t,x,a) +\mathcal{L}v(t,x,a)  + f(t,x,a,v(t,x,a),v(t,\cdot,0)-v(t,x,a))=0  \\
v(T,x,a)= g(x,a),
\end{array} \right.
\end{align}
where $D$ denotes the formal directional derivative operator
\begin{equation}\label{directional_der}
(Dv)(t,x,a):= \lim_{h \downarrow 0}\frac{v(t+h,x,a+h)-v(t,x,a)}{h}.
\end{equation}
In other  words, the presence of the directional derivative operator  \eqref{directional_der}
allows us to understand the  nonlinear Kolmogorov equation \eqref{Kolmogorov_diff} in a classical sense.
In particular, the first equality  in \eqref{Kolmogorov_diff} is understood to hold almost everywhere on $[0,\,T]$ outside of a $dt$-null set of points which can depend on $(x,a)$.

Under appropriate boundedness assumptions we have the following result:
\begin{lemma}\label{lem_existence_uniqueness_Kolmogorov_boundedness}
Suppose that $f$ and $g$ verify Hypothesis \ref{H_1} and that \eqref{g_integrability} holds; suppose, in addition, that
\begin{equation}\label{bound_f_g}
\sup_{t \in [0,\,T],\,(x,a) \in S}\Big(  \abs{g(x,a)} + \abs{f(t,x,a,0,0)} \Big) < \infty.
\end{equation}
Then the nonlinear Kolmogorov equation \eqref{Kolmogorov_int} has a unique solution $v$ in the class of measurable bounded functions.
\end{lemma}

\proof
The result follows as usual from a fixed-point argument, that we only sketch.
Let us define a map $\Gamma$ setting $v=\Gamma(w)$ where
\begin{align*}\label{eq_v_w}
& v(t,x,a) = g(x,a+T-t) + \int_{t}^{T}\mathcal{L}w(s,x,a+s-t)\,ds \\
&\qquad \qquad  + \int_{t}^{T} f(s,x,a+s-t,w(s,x,a+s-t),w(s,\cdot,0)-w(s,x,a+s-t))\,ds. \nonumber
\end{align*}
Using the Lipshitz character of $f$ and  Hypothesis \ref{hp_dati}-ii), one can show that, for some $\beta >0$ sufficiently large, the above map is a contraction in the space of bounded measurable real functions on $[0,\,T]\times S$ endowed with the supremum norm:
\begin{displaymath}
||v||_{\ast}:= \sup_{0\leqslant t\leqslant T}\sup_{(x,a)\in S}e^{-\beta(T-t)}\abs{v(t,x,a)}.
\end{displaymath}
The  unique fixed point of $\Gamma$ gives the required solution.
\endproof
Our goal is now to remove the boundedness assumption \eqref{bound_f_g}.
To this end we need to define a formula of It$\hat{\mbox{o}}$ type for the composition of the process $(X_s,\,a_s)_{s \geqslant 0}$ with functions $v$ smooth enough defined on $[0,\,T]\times S$.
Taking into account the particular form of \eqref{Kolmogorov_int}, and the fact that the second component of the process $(X_s,\,a_s)_{s \geqslant 0}$ is linear in $s$, the idea is to use in this formula the directional derivative operator $D$ given by \eqref{directional_der}.
\begin{lemma}[A formula of It$\hat{\mbox{o}}$ type]\label{Ito formula}
Let consider functions $v: [0,\,T]\times S \rightarrow \R$ such that
\begin{itemize}
\item[(i)] $\forall\, x \in K$, $\forall\, c \in [-T,\,+\infty)$, the map $t \mapsto v(t,x,t+c)$ is absolutely continuous on $[0,\,T]$, with directional derivative $D$ given by \eqref{directional_der};
\item[(ii)] for fixed $t \in [0,\,T]$, $\{ v(t+s,y,0)- v(t+s,X_{s-},a_{s-}),\, s \in [0,\,T-t],\ y \in K\}$ belongs to $\mathcal{L}^{1}_{loc}(p)$.
\end{itemize}
Then  $\P^{x,a}$-a.s., for every $t \in [0,\,T]$,
\begin{align}\label{Ito's_formula}
v(T,X_{T-t},a_{T-t}) - v(t,x,a) &=\int_{0}^{T-t}Dv(t+s,X_{s},a_{s})\,ds
+ \int_{0}^{T-t}\mathcal{L}v(t+s,X_{s},a_{s})\,ds \nonumber\\
&+ \int_{0}^{T-t}\int_{K}\left(v(t+s,y,0) - v(t+s,X_{s-},a_{s-})\right)\,q(ds,dy),
\end{align}
where the stochastic integral is a local martingale.
\end{lemma}
\proof
We proceed by reasoning as in the proof of Theorem 26.14 in \cite{Da-bo}.
We consider a function $v: [0,\,T]\times S \rightarrow \R$ satisfying (i) and (ii),
and we denote by $N_t$ the number of jumps in the interval $[0,\,t]$:
\[
N_t = \sum_{n \geqslant 1}\one_{\{T_n \leqslant t\}}.
\]
We have
\begin{align*}
v(T,X_T,a_T)-v(0,x,a) &= v(T,X_T,a_T) - v(T_{N_T},X_{T_{N_T}}, a_{T_{N_T}})+\, \sum_{n=2}^{N_T} \left\{v(T_n,X_{T_n}, a_{T_n}) -  v(T_{n-1},X_{T_{n-1}}, a_{T_{n-1}})\right\}\nonumber\\
& +\, v(T_1,X_{T_1},a_{T_1}) - v(0,x,a).
\end{align*}
Noticing that  $X_{T_{n-}}=X_{T_{n-1}}$ for all $n \in [1,\,N_T]$, $X_{T}= X_{T_{N_T}}$, and that $a_{T_{n}}=0$ for all $n \in [1,\,N_T]$, $a_{T_{1-}}=a+T_1$, and  $a_{T_{n-}}=T_n - T_{n-1}$ for all $n \in [2,\,N_T]$, we have
$$
v(T,X_T,a_T)-v(0,x,a) = I + II + III,
$$
where
\begin{align*}
I
&= (v(T_1,X_{T_1},0) - v(T_1,X_{T_1-},a_{T_1-})) + (v(T_1,x,a + T_1) - v(0,x,a))=: I' + I'',\\
II
&= \sum_{n = 2}^{N_T}
(v(T_{n},X_{T_{n}},0)  - v(T_{n},X_{T_{n}-},a_{T_{n}-}) +
 +
\sum_{n = 2}^{N_T}(v(T_{n},X_{T_{n-1}},T_{n}- T_{n-1})  - v(T_{n-1},X_{T_{n-1}},0)))=: II'+II'',\\
III
&=
v(T,X_{T},T - T_N) - v(T_N,X_{T_N},0).
\end{align*}
Let $H$ denote the $\mathcal{P}\otimes \mathcal{K}$-measurable process
\begin{displaymath}
H_s(y)= v(s,y,0)-v(s,X_{s-},a_{s-}),
\end{displaymath}
with the convention $X_{0-}=X_0$, $a_{0-}=a_0$. We have
\begin{align*}
I'+II' &= \sum_{n \geqslant 1: T_n \leqslant T} (v(T_n,X_{T_n},0)-v(T_n,X_{T_{n-}},a_{T_{n-}}))
= \sum_{n \geqslant 1: T_n \leqslant T} H_{T_n}(X_{T_N})= \int_{0}^{T}\int_{K}H_s(y)\,p(ds,dy).
\end{align*}
On the other hand, since $v$ satisfies (i) and recalling the definition \ref{directional_der} of the directional derivative operator $D$,
\begin{align*}
&I''+II''+ III
= \int_0^{T_1}\lim_{h \rightarrow 0}\frac{v(0 + h s,x,a + h s) - v(0,x,a)}{h}\,ds\\
&+\sum_{n \geqslant 2: T_n \leqslant T}\int_{T_{n-1}}^{T_{n}}\lim_{h \rightarrow 0}\frac{v(T_{n-1} + h (s-T_{n-1}),X_{T_{n-1}},a_{T_{n-1}} + h (s- T_{n-1})) - v(T_{n-1},X_{T_{n-1}},a_{T_{n-1}})}{h}\,ds\\
&+\int_{T_{N_T}}^{T}\lim_{h \rightarrow 0}\frac{v(T_{N_T} + h (s-T_{N_T}),X_{T_{N_T}},a_{T_{N_T}} + h (s- T_{N_T})) - v(T_{N_T},X_{T_{N_T}},a_{T_{N_T}})}{h}\,ds\\
& = \int_0^{T}Dv(s,X_s,a_s)\,ds.
\end{align*}
Then $\P^{x,a}$-a.s.,
\begin{align*}
&v(T,X_{T},a_{T}) - v(0,x,a)  = \int_{0}^{T} D v(s, X_s,a_s)\, ds + \int_{0}^{T}\int_{K}\left(v(s,y,0) - v(s,X_{s-},a_{s-})\right)\,p(ds,dy)\\
&\qquad = \int_{0}^{T} D v(s, X_s,a_s)\, ds +  \int_{0}^{T}\mathcal{L}v(s,X_{s},a_{s})\,ds  + \int_{0}^{T}\int_{K}\left(v(s,y,0) - v(s,X_{s-},a_{s-})\right)\,q(ds,dy),
\end{align*}
where the second equality is obtained using the identity $q(dt \,dy)= p(dt \,dy)-\lambda(X_{t-},a_{t-})\,\bar{q}(X_{t-},a_{t-},dy)\,dt$ together with the definition \eqref{L_operator} of the operator $\mathcal{L}$.

Finally, applying a shift in time, i.e. considering for every $t \in [0,\,T]$ the differential of the process $v(s+t,X_{s-},a_{s-})$ with respect to $s \in [0,\,T-t]$, the previous formula becomes: $\P^{x,a}$-a.s., for every $t \in [0,\,T]$,
\begin{align*}
v(T-t,X_{T},a_{T}) - v(t,x,a)
&= \int_{0}^{T-t} D v(s+t, X_s,a_s)\, ds +  \int_{0}^{T-t}\mathcal{L}v(s+t,X_{s},a_{s})\,ds \\
& + \int_{0}^{T-t}\int_{K}\left(v(s+t,y,0) - v(s+t,X_{s-},a_{s-})\right)\,q(ds,dy),
\end{align*}
where the stochastic integral is a local martingale thanks to condition (ii).
\endproof
We will call \eqref{Ito's_formula} the It$\hat{\mbox{o}}$ formula for $v(t+s, \cdot,\cdot) \circ {(X_{s},a_{s})}_{s \in [0,\,T-t]}$.
In differential notation:
\begin{align*}
dv(t+s,X_{s-},a_{s-}) &= Dv(t+s,X_{s-},a_{s-})\,ds \,+\,
\mathcal{L}v(t+s,X_{s-},a_{s-})\,ds \\
&\quad + \int_{K}\left(v(t+s,y,0) - v(t+s,X_{s-},a_{s-})\right)\,q(ds,dy).
\end{align*}
\begin{remark}
With respect to the classical It$\hat{\mbox{o}}$ formula, we underline that in \eqref{Ito's_formula} we have
\begin{itemize}
\item[-] the directional derivative operator $D$ instead of the usual time derivative;
\item[-] the temporal translation in the first component of $v$, i.e. we consider the differential of the process \\$v(t+s,X_{s-},a_{s-})$ with respect to $s \in [0,\,T-t]$.
        Indeed, the time-homogeneous Markov process ${(X_{s},a_{s})}_{s \geqslant 0}$ satisfies
        \begin{equation*}
        \P^{x,a}(X_{0}=x,\,a_{0}=a) = 1,
        \end{equation*}
        and the temporal translation in the first component allows us to consider $dv(t,X_t,a_t)$ for all $t \in [0,\,T]$.
\end {itemize}
\end{remark}
We go back to consider the Kolmogorov equation \eqref{Kolmogorov_int} in a more general setting. More precisely, on the functions $f$, $g$ we will only ask that they satisfy Hypothesis \ref{H_1} for every $(x,a)\in S$ and that \eqref{g_integrability} holds.
\begin{definition}\label{def_v_sol_Kolmogorov}
We say that a measurable function $v: [0,\,T] \times S \rightarrow \R$ is a solution of the nonlinear Kolmogorov equation \eqref{Kolmogorov_int}, if, for every fixed $t \in [0,\,T]$, $(x,a) \in S$,
\begin{itemize}
\item[1.] $\sperxa{\int_{0}^{T-t}\int_{K}\abs{v(t+s,y,0)-v(t+s,X_{s},a_{s})}^{2}\lambda(X_{s},a_{s})\,\bar{q}(X_{s},a_{s},dy)\,ds}<\infty$;
\item[2.] $\sperxa{\int_{0}^{T-t}\abs{v(t+s,X_{s},a_{s})}^{2}ds}<\infty$;
\item[3.] \eqref{Kolmogorov_int} is satisfied.
\end{itemize}
\end{definition}

\begin{remark}
Condition 1. is equivalent to the fact that $v(t+s,y,0) - v(t+s,X_{s-},a_{s-})$ belongs to $\mathcal{L}^{2}(p)$. Conditions 1. and 2. together are equivalent to the fact that the pair \\$\{ v(t+s,X_{s},a_{s}) , \, v(t+s,y,0) - v(t+s,X_{s-},a_{s-});\,s \in [0,\,T-t],\, y \in K \}$ belongs to the space $\mathbb{M}^{x,a}$; in particular they hold true for every measurable bounded function $v$.
\end{remark}

\begin{remark}
We need to verify the well-posedness of equation \eqref{Kolmogorov_int} for a function $v$ satisfying the condition 1. and 2. above. 
We start by noticing that, for every $(x,a)\in S$, 
$\mathbb{P}^{x,a}$-a.s.,
\begin{displaymath}
\int_{0}^{T}\int_{K}\abs{v(s,y,0)-v(s,X_{s},a_{s})}^{2}\lambda(X_{s},a_{s})\,\bar{q}(X_{s},a_{s},dy)\,ds +
\int_{0}^{T}\abs{v(s,X_{s},a_{s})}^{2}ds <\infty.
\end{displaymath}
By the law \eqref{jumpkerneldue} of the first jump  it follows that the set $\{\omega \in \Omega:\, T_{1}(\omega) >T\}$ has positive $\mathbb{P}^{x,a}$ probability, and on this set we have $X_{s-}(\omega) = x$, $a_{s-}(\omega) = a+s$. Taking such an $\omega$ we get
\begin{displaymath}
\int_{0}^{T}\int_{K}\abs{v(s,y,0)-v(s,x,a+s)}^{2}\,\lambda(x,a+s)\,\bar{q}(x,a+s,dy)\,ds +
\int_{0}^{T}\abs{v(s,x,a+s)}^{2}ds <\infty,\,\,\forall (x,a) \in S.
\end{displaymath}
Since  $\sup_{(x,a)\in S}\lambda(x,a)\bar{q}(x,a,K) < \infty$ by assumption,  H\"{o}lder's inequality implies that
\begin{align*}
\int_{0}^{T}\abs{\mathcal{L}(v(s,x,a+s))}\,ds &\leqslant
\int_{0}^{T}\int_{K}\abs{v(s,y,0)-v(s,x,a+s)}\,\lambda(x,a+s)\,\bar{q}(x,a+s,dy)\,ds \\
&\leqslant
c\left(\int_{0}^{T}\int_{K}\abs{v(s,y,0)-v(s,x,a+s)}^{2}\,\lambda(x,a+s)\,\bar{q}(x,a+s,dy)\,ds\right)^{1/2} < \infty
\end{align*}
for some constant $c$ and for all $(x,a)\in S$.
Similarly, since $\sperxa{\int_{0}^{T}\abs{f(s,X_{s},a_{s},0,0)}^2 ds}< \infty$ and arguing again on the jump time $T_{1}$, we deduce that
\begin{displaymath}
\int_{0}^{T}\abs{f(s,x,a+s,0,0)}^2\, ds< \infty,\,\,\forall (x,a)\in S;
\end{displaymath}
finally, from the Lipschitz conditions on $f$ we can conclude that
\begin{align*}
&\int_{0}^{T}\abs{f(s,x,a+s,v(s,x,a+s),v(s,\cdot,0)-v(s,x,a+s))}\,ds \\
& \leqslant
c_{1}\left(\int_{0}^{T}\abs{f(s,x,a+s,0,0)}^{2}ds\right)^{1/2} +
c_{2}\left (\int_{0}^{T}\abs{v(s,x,a+s)}^{2}ds\right)^{1/2}\\
&  +
c_{3}\left(\int_{0}^{T}\int_{K}\abs{v(s,y,0)-v(s,x,a+s)}^{2}\,\lambda(x,a+s)\,\bar{q}(x,a+s,dy)\,ds\right)^{1/2} < \infty
\end{align*}
for some constants $c_{i}$, $i = 1,2,3$, and for all $(x,a) \in S$.
Therefore, all terms occurring in equation \eqref{Kolmogorov_int} are well defined.
\end{remark}

\medskip

For every fixed $t \in [0,\,T]$ and $(x,a) \in S$, we consider now  a BSDE of the form
\begin{equation}\label{BSDE_txa}
Y^{x,a}_{s,t} + \int_{s}^{T-t}\int_{K}Z^{x,a}_{r,t} (y)\,q(dr\,dy) =
g(X_{T-t},a_{T-t}) + \int_{s}^{T-t}\,f\Big(t+r,X_{r-},a_{r-},Y^{x,a}_{r,t} ,Z^{x,a}_{r,t} (\cdot)\Big)\,dr,\,\,s \in [0,\,T-t].
\end{equation}
Then there exists a unique solution $(Y^{x,a}_{s,t} ,Z^{x,a}_{s,t} (\cdot))_{s\in [0,\,T-t]}$,  in the sense of Theorem \ref{thm: uniqueness_existence_BSDE}, and $Y^{x,a}_{0,t}$ is deterministic.
We are ready to state the main result of this section.
\begin{theorem}\label{thm_kolm}
Suppose that $f$, $g$ satisfy Hypothesis \ref{H_1} for every $(x,a) \in S$  and that \eqref{g_integrability} holds. Then for every  $t  \in [0,\,T]$, the nonlinear Kolmogorov equation \eqref{Kolmogorov_int} has a unique solution $v(t,x,a)$ in the sense of Definition \ref{def_v_sol_Kolmogorov}.

Moreover, for every fixed $t \in [0,\,T]$, for every $(x,a) \in S$ and $s \in [0,\,T-t]$  we have
\begin{align}
Y^{x,a}_{s,t} &= v(t+s,X_{s-},a_{s-}),\label{ident_1}\\
Z^{x,a}_{s,t}(y) &= v(t+s,y,0)-v(t+s,X_{s-},a_{s-})\label{ident_2},
\end{align}
so that in particular $v(t,x,a) = Y_{0,t}^{x,a}$.
\end{theorem}

\begin{remark}\label{rem_interpret_identities}
The equalities \eqref{ident_1} and \eqref{ident_2} are understood as follows.
\begin{itemize}
\item $\mathbb{P}^{x,a}$-a.s., equality \eqref{ident_1} holds for all $s \in [0,\,T-t]$.
    The trajectories of $(X_s)_{s \in [0,\,T-t]}$ are piecewise constant and càdlàg, while the trajectories of $(a_s)_{s \in [0,\,T-t]}$ are piecewise linear in $s$ (with unitary slope) and càdlàg; moreover the processes $(X_s)_{s \in [0,\,T-t]}$ and $(a_s)_{s \in [0,\,T-t]}$ have the same jump times $(T_{n})_{n\geqslant 1}$.
    Then the equality \eqref{ident_1} is equivalent to the condition
    \begin{displaymath}
    \sperxa{\int_{0}^{T-t}\abs{Y^{x,a}_{s,t}-v(t+s,X_{s},a_{s})}^{2}ds}=0.
    \end{displaymath}
\item The equality  \eqref{ident_2} holds for all $(\omega,s,y)$ with respect to the measure \\ $\lambda(X_{s-}(\omega),a_{s-}(\omega))\,\bar{q}(X_{s-}(\omega),a_{s-}(\omega),dy)\,\mathbb{P}^{x,a}(d\omega)ds$, i.e.,
    \begin{displaymath}
    \sperxa{\int_{0}^{T-t}\int_{K}\abs{Z^{x,a}_{s,t}(y)-v(t+s,y,0)+v(t+s,X_{s},a_{s})}^{2}\lambda(X_{s},a_{s})\,\bar{q}(X_{s},a_{s},dy)\,ds}=0.
    \end{displaymath}
\end{itemize}
\end{remark}

\proof
\emph{Uniqueness.} Let $v$ be a solution of the nonlinear Kolmogorov equation \eqref{Kolmogorov_int}. It follows from equality \eqref{Kolmogorov_int} itself that for every $x \in K$ and every $\tau \in [-T ,\, + \infty)$,\,\,$t \mapsto v(t,x,t+\tau)$ is absolutely continuous on $[0,\,T]$.
Indeed, applying in \eqref{Kolmogorov_int} the change of variable $\tau:= a-t$,  we obtain $\forall t \in [0,\,T]$, $\forall \tau \in [-T ,\, + \infty)$,
$$
v(t,x,t+\tau) = g(x,T+\tau) +\int_{t}^{T}\mathcal{L}v(s,x,s+\tau)\,ds+\int_{t}^{T}\, f(s,x,s+\tau,v(s,x,s+\tau),v(s,\cdot,0)-v(s,x,s+\tau))\,ds. \nonumber
$$
Then, since by assumption the process $v(t+s,y,0)-v(t+s,X_{s-},a_{s-})$ belongs to $\mathcal{L}^2(p)$, we are in a position to apply the It$\hat{\mbox{o}}$ formula \eqref{Ito's_formula} to the process $v(t+s,X_{s-},a_{s-})$, $s \in [0,\,T-t]$. We get: $\mathbb{P}^{x,a}$-a.s.,
\begin{align*}
v(t+s,X_{s-},a_{s-}) &= v(t,x,a) + \int_{0}^{s}Dv(t+r,X_{r},a_{r})\,dr
+ \int_{0}^{s}\mathcal{L}v(t+r,X_{r},a_{r})\,dr \nonumber\\
&+ \int_{0}^{s}\int_{K}\left(v(t+r,y,0) - v(t+r,X_{r},a_{r})\right)q(dr,dy),\qquad  s \in [0,\,T-t].
\end{align*}
We know that $v$ satisfies \eqref{Kolmogorov_diff}; moreover the process $X$  has piecewise constant trajectories, the process $a$ has linear trajectories in $s$,  and they have the same time jumps.
Then, $\P^{x,a}$-a.s.,
$$
Dv(t+s,X_{s-},a_{s-}) + \mathcal{L}v(t+s,X_{s-},a_{s-}) +f(t+s, X_{s-},a_{s-},v(t+s,X_{s-},a_{s-}),v(t+s,\cdot,0) - v(t+s,X_{s-},a_{s-})) = 0,
$$
for almost $s \in [0,\,T-t]$. In particular, $\P^{x,a}$-a.s.,
\begin{align*}
&v(t+s,X_{s-},a_{s-}) = v(t,x,a) + \int_{0}^{s}\int_{K}\left(v(t+r,y,0) - v(t+r,X_{r-},a_{r-})\right)q(dr,dy) \\
&\qquad - \int_{0}^{s}f(t+r, X_{r},a_{r},v(t+s,X_{s},a_{s}),v(t+r,y,0) - v(t+r,X_{r},a_{r}))\,dr, \qquad  s \in [0,\,T-t].
\end{align*}
Since $v(T,x,a)=g(x,a)$ for all $(x,a)\in S$, by simple computations we can prove that, $\forall s \in [0,\,T-t]$,
\begin{align*}
&v(t+s,X_{s-},a_{s-}) + \int_{s}^{T-t}\int_{K}\left(v(t+r,y,0) - v(t+r,X_{r-},a_{r-})\right)q(dr,dy)\\
& \qquad  = g(X_{T-t},a_{T-t}) + \int_{s}^{T-t}\,f(t+r, X_{r},a_{r},v(t+r,X_{r},a_{r}),v(t+r,y,0) - v(t+r,X_{r},a_{r}))\,dr \nonumber.
\end{align*}
Since
the pairs $(Y^{x,a}_{s,t},Z^{x,a}_{s,t}(\cdot))_{s \in [0,\,T-t]}$ and $(v(t+s,X_{s-},a_{s-}),v(t+s,y,0)-v(t+s,X_{s-},a_{s-}))_{s \in [0,\,T-t]}$ are both solutions to the same BSDE under $\P^{x,a}$, they coincide as members of the space $\mathbb{M}^{x,a}$. It follows that equalities \eqref{ident_1} and \eqref{ident_2} hold. In particular,  $v(t,x,a)= Y_{0,t}^{x,a}$, and this yields the uniqueness of the solution.

\medskip

\emph{Existence.}
We proceed by   an approximation argument, following the same lines of the proof of Theorem 4.4 in \cite{CoFu-m}.
We recall that, by Theorem \ref{thm: uniqueness_existence_BSDE}, for every fixed $t \in [0,\,T]$, the BSDE \eqref{BSDE_txa} has a unique solution $(Y^{x,a}_{s,t},Z^{x,a}_{s,t}(\cdot))_{s \in [0,\,T-t]}$ for every $(x,a) \in S$; moreover, $Y_{0,t}^{x,a}$ is deterministic, i.e., there exists a real number, denoted by $v(t,x,a)$, such that $\P^{x,a}(Y_{0,t}^{x,a}= v(t,x,a))=1$.
At this point, we set  $f^{n} = (f \wedge n)\vee (-n)$ and  $g^{n} = (g \wedge n)\vee (-n)$ as the truncations of $f$ and $g$ at level $n$. By Lemma \ref{lem_existence_uniqueness_Kolmogorov_boundedness},  for $t \in [0,\,T]$, $(x,a)\in S$, equation
\begin{align}\label{Kolmogorov_truncated}
 v^{n}(t,x,a) &= g^{n}(x,a+T-t) + \int_{t}^{T}\mathcal{L}v^{n}(s,x,a+s-t)\,ds \nonumber\\
& \quad  + \int_{t}^{T}\, f^{n}(s,x,a+s-t,v^{n}(s,x,a+s-t),v^{n}(s,\cdot,0)-v^{n}(s,x,a+s-t))\,ds.
\end{align}
admits a unique bounded measurable solution $v^{n}$.
In particular,  the first part of the proof yield the following identifications:
\begin{align*}
v^{n}(t,x,a) &= Y^{x,a,n}_{0,t},\\
v^{n}(t+s,X_{s-},a_{s-}) &= Y^{x,a,n}_{s,t},\\
v^{n}(t+s,y,0) - v^{n}(t+s,X_{s-},a_{s-}) &= Z^{x,a,n}_{s,t}(y),
\end{align*}
in the sense of Remark \ref{rem_interpret_identities}, where $(Y_{s,t}^{x,a,n},Z_{s,t}^{x,a,n}(\cdot))_{s \in [0,\,T-t]}$ is the unique solution to the BSDE
$$
Y^{x,a,n}_{s,t} + \int_{s}^{T-t}\int_{K}Z^{x,a,n}_{r,t}(y)\,q(dr\,dy)= g^{n}(X_{T-t},a_{T-t}) + \int_{s}^{T-t}\,f^{n}\left(t+r,X_{r},a_{r},Y^{x,a,n}_{r,t},Z^{x,a,n}_{r,t}(\cdot)\right)\,dr,
$$
for all $s\in [0,\,T-t]$.
Recalling \eqref{BSDE_txa} and applying Theorem \ref{thm: uniqueness_existence_BSDE}, we deduce that,  for  some constant $c$,
\begin{align}\label{convergence}
&\sup_{s \in [0,\,T-t]}\E^{x,a}\left[|Y^{x,a}_{s,t}-Y^{x,a,n}_{s,t}|^2\right] +
\E^{x,a}\left[\int_{0}^{T-t}|Y^{x,a}_{s,t}-Y^{x,a,n}_{s,t}|^2ds\right] \nonumber \\
&\qquad  + \E^{x,a}\left[\int_{0}^{T-t}\int_{K}|Z^{x,a}_{s,t}(y)-Z^{x,a,n}_{s,t}(y)|^2\lambda(X_{s},a_{s})\,\bar{q}(X_{s},a_{s},dy)\,ds\right]\nonumber\\
&
\leqslant c \E^{x,a}\left[|g(X_{T-t},a_{T-t})-g^{n}(X_{T-t},a_{T-t})|^2\right] \nonumber \\
&\qquad+c\E^{x,a}\left[\int_{0}^{T-t}|f(t+s,X_{s},a_{s},Y^{x,a}_{s,t},Z^{x,a}_{s,t}(\cdot))-f^{n}(t+s,X_{s},a_{s},Y^{x,a}_{s,t},Z^{x,a}_{s,t}(\cdot))|^2ds\right] \longrightarrow 0,
\end{align}
where the two final terms tend to zero by monotone convergence.
In particular \eqref{convergence} yields
\begin{displaymath}
|v(t,x,a)-v^{n}(t,x,a)|^{2} = |Y^{x,a}_{0,t}-Y^{x,a,n}_{0,t}|^2 \leqslant
\sup_{s \in [0,\,T-t]}\E^{x,a}\left[|Y^{x,a}_{s,t}-Y^{x,a,n}_{s,t}|^2\right]\longrightarrow 0,
\end{displaymath}
and therefore $v$ is a measurable function.
At this point, applying the Fatou Lemma we get
\begin{align*}
&\sperxa{\int_{0}^{T-t}\abs{Y^{x,a}_{s,t}-v(t+s,X_{s},a_{s})}^{2}\,ds} \\
&+
\sperxa{\int_{0}^{T-t}\int_{K}\abs{Z^{x,a}_{s,t}(y)-v(t+s,y,0)+v(t+s,X_{s},a_{s})}^{2}\,\lambda(X_{s},a_{s})\,\bar{q}(X_{s},a_{s},dy)\,ds}\\
&\leqslant
\liminf_{n \rightarrow \infty}\sperxa{\int_{0}^{T-t}\abs{Y^{x,a}_{s,t}-v^{n}(t+s,X_{s},a_{s})}^{2}\,ds}\\
&+
\liminf_{n \rightarrow \infty}\sperxa{\int_{0}^{T-t}\int_{K}\abs{Z^{x,a}_{s,t}(y)-v^n(t+s,y,0)+v^n(t+s,X_{s},a_{s})}^{2}\,\lambda(X_{s},a_{s})\,\bar{q}(X_{s},a_{s},dy)\,ds}\\
& = \liminf_{n \rightarrow \infty}\sperxa{\int_{0}^{T-t}\abs{Y^{x,a}_{s,t}-Y_{s,t}^{x,a,n}}^{2}\,ds} \\
&+
\liminf_{n \rightarrow \infty}\sperxa{\int_{0}^{T-t}\int_{K}\abs{Z^{x,a}_{s,t}(y)-Z_{s,t}^{x,a,n}(y)}^{2}\lambda(X_{s},a_{s})\,\bar{q}(X_{s},a_{s},dy)\,ds}= 0
\end{align*}
by \eqref{convergence}. The above calculations show  that \eqref{ident_1} and \eqref{ident_2} hold. Moreover, they  imply that
\begin{align*}
&\sperxa{\int_{0}^{T-t}\abs{v(t+s,X_{s},a_{s})}^{2}ds}  +
\sperxa{\int_{0}^{T-t}\int_{K}\abs{v(t+s,y,0)-v(t+s,X_{s},a_{s})}^{2}\lambda(X_{s},a_{s})\,\bar{q}(X_{s},a_{s},dy)\,ds}\\
& =
\sperxa{\int_{0}^{T-t}\abs{Y^{x,a}_{s,t}}^{2}ds}+
\sperxa{\int_{0}^{T-t}\int_{K}\abs{Z^{x,a}_{s,t}(y)}^{2}\lambda(X_{s-},a_{s-})\,\bar{q}(X_{s},a_{s},dy)\,ds}< \infty,
\end{align*}
that accords to requirement of Definition \ref{def_v_sol_Kolmogorov}.

It  remains to show that $v$ satisfies \eqref{Kolmogorov_int}. This would follow from a passage to the limit in \eqref{Kolmogorov_truncated}, provided we  show that
\begin{equation}\label{convergence_Kolm_trunc_1}
\int_{t}^{T}\mathcal{L}v^{n}(s,x,a+s-t)ds  \rightarrow \int_{t}^{T}\mathcal{L}v(s,x,a+s-t)ds,
\end{equation}
and
\begin{align}\label{convergence_Kolm_trunc_2}
&\int_{t}^{T} f^{n}(s,x,a+s-t,v^{n}(s,x,a+s-t),v^{n}(s,\cdot,0)-v^{n}(s,x,a+s-t))\,ds \nonumber\\
& \qquad \qquad \rightarrow \int_{t}^{T} f(s,x,a+s-t,v(s,x,a+s-t),v(s,\cdot,0)-v(s,x,a+s-t))\,ds.
\end{align}
To prove \eqref{convergence_Kolm_trunc_1}, we observe that
\begin{align*}
&\E^{x,a}\abs{\int_{0}^{T-t}\mathcal{L}v(t+s,X_{s-},a_{s-})\,ds - \int_{0}^{T-t}\mathcal{L}v^{n}(t+s,X_{s-},a_{s-})\,ds}
\\
&=
\E^{x,a}\abs{\int_{0}^{T-t}\int_{K}(Z^{x,a}_{s,t}-Z_{s,t}^{x,a,n})\,\lambda(X_{s},a_{s})\,\bar{q}(X_{s},a_{s},dy)\,ds}
\\
& \leqslant (T-t)^{1/2}\sup_{x,a}[\lambda(x,a)\,\bar{q}(x,a,K)]^{1/2}\left (\sperxa{\int_{0}^{T-t}\int_{K}\abs{Z^{x,a}_{s,t}-Z_{s,t}^{x,a,n}}\,\lambda(X_{s},a_{s})\,\bar{q}(X_{s},a_{s},dy)\,ds}\right)^{1/2}\rightarrow 0,
\end{align*}
by \eqref{convergence}.
Then,  for a subsequence (still denoted $v^n$) we get 
\[
\int_{0}^{T-t}\mathcal{L}v^n(t+s,X_{s},a_{s})\,ds \rightarrow \int_{0}^{T-t}\mathcal{L}v(t+s,X_{s},a_{s})\,ds, \quad \P^{x,a}\textup{-a.s}.
\]
Recalling the law \eqref{jumpkerneldue} of the first jump $T_1$, we see that the set $\{\omega \in \Omega:\, T_{1}(\omega) >T\}$ has positive $\mathbb{P}^{x,a}$ probability, and on this set we have $X_{s-}(\omega) = x$, $a_{s-}(\omega) = a+s$. Choosing such an $\omega$ we have
\begin{displaymath}
\int_{0}^{T-t}\mathcal{L}v^n(t+s,x,a+s)ds \rightarrow \int_{0}^{T-t}\mathcal{L}v(t+s,x,a+s)ds,
\end{displaymath}
i.e., by a translation of $t$ in the temporal line,
\begin{displaymath}
\int_{t}^{T}\mathcal{L}v^n(s,x,a+s-t)ds \rightarrow \int_{t}^{T}\mathcal{L}v(s,x,a+s-t)ds.
\end{displaymath}

To show \eqref{convergence_Kolm_trunc_2}, we compute
\begin{align*}
&\E^{x,a}\abs{\int_{0}^{T-t}f(t+s,X_{s},a_{s},Y^{x,a}_{s,t},Z^{x,a}_{s,t})-f^n(t+s,X_{s},a_{s},Y_{s,t}^{x,a,n},Z_{s,t}^{x,a,n}))\,ds}
\\
&\leqslant \sperxa{\int_{0}^{T-t}\abs{f(t+s,X_{s},a_{s},Y^{x,a}_{s,t},Z^{x,a}_{s,t})-f^n(t+s,X_{s},a_{s},Y^{x,a}_{s,t},Z^{x,a}_{s,t})}\,ds}\\
& + \sperxa{\int_{0}^{T-t}\abs{f^{n}(t+s,X_{s},a_{s},Y^{x,a}_{s,t},Z^{x,a}_{s,t})-f^n(t+s,X_{s},a_{s},Y_{s,t}^{x,a,n},Z_{s,t}^{x,a,n})}\,ds}.
\end{align*}
The first integral term in the right-hand side tends to zero by monotone convergence. At this point, we notice that    $f^n$ is a truncation of $f$, and therefore it satisfies the Lipschitz condition \eqref{f_inequality} with the same constants $L$, $L'$, independent of $n$. This yields the following estimate  for  the second integral:
\begin{align*}
& L' \,\sperxa{\int_{0}^{T-t}\abs{Y^{x,a}_{s,t}-Y_{s,t}^{x,a,n}}ds}+ L\,\sperxa{\int_{0}^{T-t}\Big(\int_{K}\abs{Z^{x,a}_{s,t}(y)-Z_{s,t}^{x,a,n}(y)}^2\lambda(X_{s},a_{s})\,\bar{q}(X_{s},a_{s},dy)\Big)^{1/2}ds}\\
& \leqslant L' \, \left( \,(T-t)\,\sperxa{\int_{0}^{T-t}\abs{Y^{x,a}_{s,t}-Y_{s,t}^{x,a,n}}^2\,ds} \right)^{1/2}\\
& + L\,\left(\,(T-t)\,\sperxa{\int_{0}^{T-t}\int_{K}\abs{Z^{x,a}_{s,t}(y)-Z_{s,t}^{x,a,n}(y)}^2\lambda(X_{s},a_{s})\,\bar{q}(X_{s},a_{s},dy)\,ds}\,\right)^{1/2},
\end{align*}
which tends to zero, again by \eqref{convergence}. Considering 
a subsequence (still denoted $v^{n}$) we get,
\begin{eqnarray*}
&&\int_{0}^{T-t}f^{n}(t+s,X_{s},a_{s},v^{n}(t+s,X_{s},a_{s}),v^{n}(t+s,y,0)-v^{n}(t+s,X_{s},a_{s}))\,ds
\\
&& \qquad \rightarrow
\int_{0}^{T-t}f(t+s,X_{s},a_{s},v(t+s,X_{s},a_{s}),v(t+s,y,0)-v(t+s,X_{s},a_{s}))\,ds,\,\, \P^{x,a}\mbox{-a.s.}
\end{eqnarray*}
Choosing also in this case an $\omega$ in the set $\{\omega \in \Omega: \, T_{1}(\omega)> T\}$, we find
\begin{eqnarray*}
&&\int_{0}^{T-t}f^{n}(t+s,x,a+s,v^{n}(t+s,x,a+s),v^{n}(t+s,y,0)-v^{n}(t+s,x,a+s))\,ds
\\
&& \qquad \rightarrow
\int_{0}^{T-t}f(t+s,x,a+s,v(t+s,x,a+s),v(t+s,y,0)-v(t+s,x,a+s))\,ds,
\end{eqnarray*}
and a change of temporal variable allows to prove that   \eqref{Kolmogorov_int} holds, and to conclude the proof.
\endproof

We finally introduce the Hamilton-Jacobi-Bellman (HJB) equation associated to the control problem considered in Section \ref{section semi_markov_HJB}:
for every $t \in [0,\,T]$ and $(x,a)\in S$,
\begin{equation}\label{HJB}
v(t,x,a) = g(x,a+T-t) + \int_{t}^{T}\mathcal{L}v(s,x,a+s-t)\,ds + \int_{t}^{T} f(s,x,a+s-t,v(s,\cdot,0)-v(s,x,a+s-t))\,ds,
\end{equation}
where $\mathcal{L}$ denotes the operator introduced in \eqref{L_operator}, $f$ is the Hamiltonian function defined by \eqref{hamilton_function} and $g$ is the terminal cost.
Since \eqref{HJB}  is a nonlinear Kolmogorov equation of the form \eqref{Kolmogorov_int}, we can apply Theorem \ref{thm_kolm} and conclude that the value function and an optimal control law can be represented by means of the HJB solution $v(t,x,a)$.
\begin{corollary}
Let  Hypotheses \ref{hyp_control} and \ref{hyp_assumed_min} hold. For every fixed $t \in [0,\,T]$, for every $(x,a) \in S$ and $s \in [0,\,T-t]$, there exists a unique solution $v$ to the HJB equation \eqref{HJB}, satisfying
\begin{eqnarray*}
v(t+s,X_{s-},a_{s-}) &=& Y^{x,a}_{s,t},\\
v(t+s,y,0)-v(t+s,X_{s-},a_{s-})&=& Z^{x,a}_{s,t}(y),
\end{eqnarray*}
where the above equalities are understood as explained in Remark \ref{rem_interpret_identities}.\\
In particular an optimal control is given by the formula
$$
u^{*\,t,x,a}_{s} \in \Gamma(t+s, X_{s-}, a_{s-}, v(t+s,\cdot,0)-v(t+s,X_{s-},a_{s-})),
$$
while the value function coincides with $v(t,x,a)$, i.e. $J(t,x,a,u^{*\,t,x,a}(\cdot)) = v(t,x,a)= Y_{0,t}^{x,a}$.
\end{corollary}


\bigskip

\footnotesize
\begin{thebibliography}{11}


\bibitem{BaBuPa} Barles G., Buckdahn R., Pardoux E.
Backward stochastic differential equations and integral-partial
differential equations. Stochastics Stochastics Rep. 60 (1997),
no. 1-2, 57-83.

\bibitem{Be}
Becherer D.  Bounded solutions to backward SDEs with jumps for utility optimization and indifference hedging. Ann. of App. Prob. 16 (2006), 2027-2054.

\bibitem{Bo-Va-Wo} Boel R., Varaiya P., Wong E.  Martingales on jump processes;
Part I: Representation results; Part II: Applications; SIAM J. Control 13,
999-1061.


\bibitem{B} Br\'emaud  P. Point processes and queues, Martingale dynamics.
Springer Series in Statistics. Springer (1981).

\bibitem{chito} Chitopekar S. S.  Continuous time Markovian sequential control processes. SIAM J. Control 7 (1969), 367-389.



\bibitem{CoFu-mpp} Confortola F., Fuhrman M. Backward stochastic differential equations and optimal control of marked point processes. SIAM J. Control Optimization 51(5) (2013), 3592-3623.

\bibitem{CoFu-m} Confortola F., Fuhrman M. Backward stochastic differential equations associated to jump Markov processes and their applications. Stochastic Processes and their Applications 124 (2014), 289-316.

\bibitem{CFJ} Confortola F., Fuhrman M., Jacod J.  Backward stochastic differential equations driven by a marked point process: an elementary approach, with an application to optimal control (2014), preprint.

\bibitem{Cre-Mat} Crépey S., Matoussi A. Reflected and doubly reflected BSDEs with jumps, Ann. of Appl. Prob. 18 (2008), 2041-2069.



\bibitem{Da-bo} Davis M.H.A. Markov models and optimization.
Monographs on Statistics and Applied Probability 49 (1993), Chapman $\&$
Hall.

\bibitem{Da-Fa} Davis M.H.A., Farid M. Piecewise deterministic processes and viscosity solutions.
McEneaney, W. M. et al. (ed) Stochastic Analysis, Control Optimization and Applications. A Volume in Honor of W. H. Fleming on Occasion of His 70th Birthday, Birkh\"auser (1999), 249-268.

\bibitem{Dem} Dempster M. A. H.
Optimal control of piecewise deterministic Markov processes. Applied stochastic analysis (London, 1989), 303-325,
Stochastics Monogr. 5, Gordon and Breach (New York 1991).

\bibitem{ElK} El Karoui N.
Les aspects probabilistes du contr\^ole stocastique. [The probabilistic aspects of stochastic control] Ninth Saint Flour Probabilistic Summer School,-1979. (Saint Flour, 1979), 73-238, Lecture Notes in Math. 876 (1981), Springer.

\bibitem{E} Elliott R.J. Stochastic Calculus and its Applications. Springer,
(1982).

\bibitem{G-S} Gihman I. I., Skorohod A. V.
The Theory of Stochastic Processes II.
Springer-Verlag Berlin Heidelberg New-York Tokyo (1983).


\bibitem{Howard} Howard R. A. Dynamic Probabilistic Systems.
John Wiley, New York (1971).


\bibitem{J} Jacod J. Multivariate point processes: predictable projection,
Radon-Nikodym derivatives, representation of martingales.
Z. Wahrscheinlichkeitstheorie und Verw. Gebiete 31 (1974/75), 235-253.


\bibitem{Jewell} Jewell W. S. Markov-renewal programming I and II.
Operational Res. 11 (1963), 938-971, 938-971.

\bibitem{KaTa-Po-Zh_1}
Kazi-Tani N.,  Possama\"i D., Zhou C.
Second Order BSDEs  with jumps, Part I: Aggregation and uniqueness preprint.

\bibitem{KaTa-Po-Zh_2}
Kazi-Tani N.,  Possama\"i D., Zhou C.
Second Order BSDEs  with jumps, Part II: Existence and applications, preprint.

\bibitem{KhMaPhZh}
Kharroubi I.,  Ma J., Pham H., Zhang J.
Backward SDEs with constrained jumps and quasi-variational inequalities.
Ann. Probab. 38(2)
(2010), 794-840.

\bibitem{LenYam}
Lenhart S.-M.,  Yamada N.
Perron's method for viscosity solutions associated with piecewise-deterministic processes.
Funkcialaj Ekvacioj  34
(1991),  173-186.

\bibitem{LimQuenez}
Lim T.,  Quenez M.-C.
Exponential utility maximization and indifference price in an incomplete market with defaults.
(2010), preprint.

\bibitem{Osaki} Osaki S., Mine H. Linear programming algorithms for semi-Markovian decision processes. J. Math. Anal. 22 (1968), pp. 356-381.


\bibitem{Re-Yor} Revuz D.,  Yor M. Continuous Martingales and Brownian Motion.
Grundlehren der mathematischen Wissenschaften. Springer, third
edition (1999).


\bibitem{Ross} Ross S. M. Applied Probability Models with Optimization Applications. Holden-Day, San Francisco (1970).

\bibitem{Roy} Royer M. Backward stochastic differential equations with jumps and
related nonlinear expectations. Stochastic Processes and their Applications, 116(10) (2006), 1358-1376.

\bibitem{St1} Stone L. D.  Necessary and sufficient conditions for optimal control of semi-Markov jump processes. SIAM J. Control Optim. 11(2) (1973), 187-201.


\bibitem{TaLi}
Tang S.,  Li X. Necessary Conditions for Optimal Control of
Systems with Random Jumps. SIAM J. Control Optim. 32 (1994),
1447-1475.


\bibitem{Ver} Vermes D. Optimal control of piecewise deterministic Markov process. Stochastics 14(3) (1985), 165-207.


\bibitem{Xia} Xia J. Backward stochastic differential equations with random measures.
Acta Mathematicae Applicatae Sinica 16(3) (2000), 225-234.

\end{thebibliography}
\end{document}